%% file: laminations_and_circles.tex
\documentclass[12pt]{amsart}
\usepackage{utopia, euler}

\usepackage{amssymb}
\usepackage{amsmath}

\usepackage{xypic}
\usepackage{graphicx}
\usepackage{psfrag}
\usepackage{longtable}

\def\RCS$#1: #2 ${\expandafter\def\csname RCS#1\endcsname{#2}}
\RCS$Revision: 1.134 $
\RCS$Date: 2002/07/21 15:23:03 $

\renewcommand{\today}{\number\year /\number\month /\number\day}

  \usepackage{fullpage}
\newcommand{\R}{{\mathbb R}}
\newcommand{\Z}{{\mathbb Z}}
\renewcommand{\H}{{\mathbb H}}
\newcommand{\bdry}{\partial}
\newcommand{\iso}{\cong}

\newcommand{\maps}{\colon\thinspace}
\newcommand{\un}{{\text{univ}}}
\newcommand{\homeo}{{\text{Homeo}}}
\newcommand{\til}{\widetilde}
\newcommand{\F}{{\mathscr F}}
\newcommand{\G}{{\mathscr G}}

\DeclareMathOperator{\vol}{vol}

\newcommand{\PSL}[2]{\mathrm{PSL}_{#1} #2}

\newcommand{\spandef}[2]{{  \left\langle  {#1}  \ \left| \   {#2} \right. \right\rangle }}
\newcommand{\setdef}[3][]{{  #1\{  {#2}  \ #1| \   {#3} #1\} }}
\newcommand{\mtext}[1]{\quad\mbox{#1}\quad}
\newcommand{\pair}[1]{\left\langle #1 \right\rangle}

\newcommand{\SnapPea}{\texttt{SnapPea}}

\newcommand{\constructP}{\texttt{constructP}}

\newtheoremstyle{plain}{}{}{\slshape}{}{\bfseries}{.}{0.5em}{}
\theoremstyle{plain} 
\swapnumbers

\newtheorem{Theorem}{Theorem}[section]

\newtheorem{Lemma}[Theorem]{Lemma}
\newtheorem{Corollary}[Theorem]{Corollary}
\newtheorem{Claim}[Theorem]{Claim}
\newtheorem{Question}[Theorem]{Question}

\newtheorem{defn}[Theorem]{Definition}

\newtheorem{filling_lemma}[Theorem]{Filling Lemma}
\newtheorem{leaf_pocket}[Theorem]{Leaf Pocket Theorem for Laminations}

\newtheorem{maintheorem}{Theorem}
\newtheorem{fillingforward}{Filling Lemma}
\newtheorem{weeksforward}{Theorem}

  \theoremstyle{remark}

\newtheorem{rmk}[Theorem]{Remark}
\newtheorem{exa}[Theorem]{Example}

\newcommand{\halmos}{}

\makeatletter
  \let\c@Theorem=\c@subsection
  \let\c@figure=\c@subsection
  \let\p@figure=\p@subsection
  
  \let\cl@figure=\cl@subsection
\makeatother

\begin{document}

\title{Laminations and groups of homeomorphisms of the circle}

\author[Calegari]{Danny Calegari}
\address{Department of Mathematics, California Institute of Technology, Pasadena, CA 91125, USA}
\email{dannyc@math.harvard.edu}

  \author[Dunfield]{Nathan M. Dunfield}
\address{Department of Mathematics \\ Harvard University \\ Cambridge MA 02138, USA}
\email{nathand@math.harvard.edu}
\date{Version: \RCSRevision, Compile: \today, Last commit:
    \RCSDate}

\begin{abstract} 
  If $M$ is an atoroidal 3-manifold with a taut foliation, Thurston
  showed that $\pi_1(M)$ acts on a circle.  Here, we show that some
  other classes of essential laminations also give rise to actions on
  circles.  In particular, we show this for tight essential
  laminations with solid torus guts.  We also show that pseudo-Anosov
  flows induce actions on circles.  In all cases, these actions can be
  made into faithful ones, so $\pi_1(M)$ is isomorphic to a subgroup
  of $\homeo (S^1)$.  In addition, we show that the fundamental group
  of the Weeks manifold has no faithful action on $S^1$.  As a
  corollary, the Weeks manifold does not admit a tight essential
  lamination, a pseudo-Anosov flow, or a taut foliation.  Finally, we
  give a proof of Thurston's universal circle theorem for taut
  foliations based on a new, purely topological, proof of the
  Leaf Pocket Theorem.
\end{abstract}

\maketitle

\setcounter{tocdepth}{1}
\tableofcontents

\section{Introduction}

Let $M$ be an atoroidal 3-manifold with a taut foliation $\F$.  The
foliation $\F$ gives rise to actions of $\pi_1(M)$ on 1-manifolds in
two ways: the action on the space of leaves in the universal cover,
and the action on a universal circle.  These actions are very useful
for understanding $\F$.  If $M$ is an atoroidal $3$-manifold with an
essential lamination, one still has an action of $\pi_1(M)$ on the
space of leaves in the universal cover, but this space is just an
order tree, not a 1-manifold.  Our main goal here is to give analogues
of the second kind of action: we construct actions of fundamental
groups of $3$-manifolds on circles arising from certain kinds of
essential laminations.  Before giving precise statements, we'll
discuss the case of taut foliations.

Let $\til{M}$ denote the universal cover of $M$, and $\til{\F}$ the
foliation of $\til{M}$ covering a taut foliation $\F$ of $M$.  First,
consider the leaf space $L$ of $\til{\F}$. The space $L$ is simply
connected and locally homeomorphic to $\R$, but is typically
non-Hausdorff.  We will still refer to $L$ as a 1-manifold.  There is
a natural action of $\pi_1(M)$ on $L$, and this action has no global
fixed point.  In fact, after changing $\F$ slightly by a monotone
equivalence, we can ensure that the homomorphism $\pi_1(M) \to
\homeo(L)$ is injective (see Theorem~\ref{fol_action_on_leafspace}).
Thus $\pi_1(M)$ is isomorphic to a group of homeomorphisms of a
1-manifold.

Second, Thurston has shown that the foliation $\F$ gives rise to
another action of $\pi_1(M)$ on a 1-manifold, namely on a
\emph{universal circle}.  This circle, $S^1_\un$, is constructed by
collating the circles at infinity of the leaves of $\til{\F}$.  As in
the previous case, the action on a universal circle is faithful
(Theorem~\ref{universal_circle_action_faithful}).  Thurston's paper
\cite{ThurstonCirclesII} on this topic is largely unwritten, so we
provide a complete proof here.  The main technical tool for
proving the existence of a universal circle is Thurston's
\emph{Leaf Pocket Theorem}.  This theorem says, very roughly, that the
leaves of $\F$ come together in many directions.  Thurston's proof of
this theorem is analytic, and uses theorems of L.~Garnett on the
existence and quality of harmonic transverse measures for foliations.
We give a new and purely topological proof of a slight variation of
the Leaf Pocket Theorem, Theorem~\ref{leaf_pocket_theorem}, which
applies more generally to essential laminations.  In
Section~\ref{universal_circles}, we use
Theorem~\ref{leaf_pocket_theorem} to prove the existence of universal
circles.

Now we will outline our results for essential laminations.  For
background and more detailed definitions see
Section~\ref{essential_laminations}.  An essential lamination
$\Lambda$ of a 3-manifold $M$ is \emph{tight} if the leaf space of
$\til{\Lambda}$ is Hausdorff.  Equivalently, leaves of $\til{\Lambda}$
are uniformly properly embedded in $\til{M}$.  The complementary
regions of $\Lambda$ can be partitioned into \emph{interstitial
  regions} and \emph{guts}.  The interstitial regions are (typically
noncompact) $I$-bundles whose $\bdry I$-bundles lie on $\Lambda$, and
the compact guts make up the rest of the complement.

Our main result is
\begin{maintheorem}
  Let $M$ be an atoroidal 3-manifold containing an essential
  lamination $\Lambda$.  If $\Lambda$ is tight and all the gut regions
  are solid tori, then $\pi_1(M)$ has a faithful action on $S^1$.
\end{maintheorem}

Here's an outline of the proof.  We can assume that $\Lambda$ has
non-empty guts (i.e.~is \emph{genuine}), as otherwise it can be blown down
to a taut foliation which gives rise to an action on a universal circle.  The
first step in the proof is the following lemma, which is of
independent interest:
\begin{fillingforward}
  Let $M$ be a $3$-manifold and $\Lambda$ a genuine lamination all of
  whose guts are solid tori.  Then there is a genuine lamination
  $\overline{\Lambda}$ containing $\Lambda$, whose complementary
  regions are all ideal polygon bundles over $S^1$.  Moreover, if
  $\Lambda$ is tight, so is $\overline{\Lambda}$.
\end{fillingforward}
A lamination such as $\overline{\Lambda}$ whose complementary regions
are all ideal polygon bundles is said to be \emph{very full}.  A good
example to keep in mind of a very full lamination is the suspension of
the stable lamination of a pseudo-Anosov surface homeomorphism in the
corresponding surface bundle.

The Filling Lemma~\ref{filling_lemma} reduces the proof of
Theorem~\ref{main_thm} to the case where $\Lambda$ is tight and very
full, which is handled in Theorem~\ref{order_tight_full}.  The basic
idea is to flatten out $\til{\Lambda}$ into a geodesic lamination
$\lambda$ of $\H^2$, taking care so that $\lambda$
inherits a group action of $\pi_1(M)$ which is equivariant with
respect to the map $\til{\Lambda} \to \lambda$.  The action of
$\pi_1(M)$ on $\lambda$ induces an action on the circle at infinity of
$\H^2$, which is the desired circle action.

We can lessen the requirement that $M$ be tight somewhat, when the
non-Haus\-dorff behavior of the lamination is not too bad.  See
Section~\ref{sec_non_tight} and Theorem~\ref{thm_ord_cataclysms}.  In
particular Theorem~\ref{thm_ord_cataclysms} applies to the stable or
unstable lamination of a pseudo-Anosov flow.  Thus
(Corollary~\ref{cor_anosov_flow}) fundamental groups of 3-manifolds with pseudo-Anosov flows
also act faithfully on the circle.

\subsection{Nonexistence results}  

Given the results above, it is natural to ask the question: Which
atoroidal 3-manifolds have faithful circle actions?  We show that not
every 3-manifold has such an action.  In particular, for the closed
hyperbolic 3-manifold with smallest known volume, we have:
\begin{weeksforward}
  The fundamental group of the Weeks manifold does not act faithfully
  on the circle.
\end{weeksforward}
As a corollary of this and the above theorems on the existence of
circle actions, it follows that the Weeks manifold does not admit a
tight essential lamination with solid torus guts, a pseudo-Anosov
flow, or a taut foliation.  Using Agol's volume estimates for
manifolds with tight laminations \cite{AgolVolume}, one can show that
any tight lamination in the Weeks manifold would have solid torus
guts.  So we have the stronger conclusion that the Weeks manifold
contains no tight essential lamination
(Corollary~\ref{weeks_nonexistence}).  We conjecture that the Weeks
manifold contains no essential lamination at all, but this remains
open.   Also, the fundamental group of the Weeks manifold appears to be
the first known example of a rank-1 lattice which cannot act
faithfully on the circle.  In contrast, Witte has shown that there are
many higher rank lattices with this property \cite{Witte94}.

Examples of hyperbolic manifolds without taut foliations were already
known from the breakthrough work of Roberts, Shareshian, and Stein
\cite{RobertsShareshianStein}.  In fact, they constructed infinitely
many such examples by looking at certain Dehn fillings of punctured
torus bundles with negative trace.  They showed the non-existence of
taut foliations by proving that the fundamental groups of these
manifolds can't act on a simply-connected 1-manifold without a global
fixed point.  Given the multitude of non-Hausdorff 1-manifolds, it is
remarkable that such a proof can be made to work.  One of our original
motivations for this paper was to reprove non-existence of taut
foliations by studying only actions on Hausdorff 1-manifolds, in
particular actions on $S^1$.   

The technique used for proving Theorem~\ref{weeks_circle} is this.
Let $G$ be the fundamental group of the Weeks manifold.  First, we
pass to a finite index subgroup $N$ where the action on $S^1$ lifts to
an action on $\R$.  A faithful action of $N$ on $\R$ is equivalent to
a left-invariant total order on $N$ (see Section~\ref{sec_orderings}).
The non-existence of such an order is then shown by considering
various possibilities for which elements of $N$ satisfy $g > 1$.  For
some $3$-manifold groups this can be done algorithmically as discussed
in Section~\ref{algorithmic_issues}.  While the Weeks manifold was the
only case where we managed to show that the fundamental group doesn't
act faithfully on the circle, we found many examples in the
Hodgson-Weeks census whose fundamental groups can't act faithfully on
$\R$ (see Section~\ref{census_exs}).

\subsection{Further Questions}

Because essential laminations are very common, these results show that
many 3-manifold groups are subgroups of $\homeo(S^1)$.  So a natural
question is: Suppose that $\pi_1(M)$ is a subgroup of $\homeo(S^1)$.
What does this tell us about the algebraic properties of $\pi_1(M)$?
Of course, a general finitely generated subgroup of $\homeo(S^1)$ can
be quite strange, e.g it can contain Thompson's infinite simple group.
However, the actions we construct have additional special properties; for
instance, they preserve the endpoints of a special kind of geodesic
lamination of $\H^2$.  So one can hope that the existence of such an
action implies something interesting about $\pi_1(M)$.

\subsection{Acknowledgments}

The first author was at Harvard University when this work was done, and
was partially supported by an NSF VIGRE grant.  The second author was
partially supported by an NSF Postdoctoral Fellowship.  We thank Bill
Thurston, Rachel Roberts, and Ian Agol for useful conversations.  We
also thank Sergio Fenley and the referee for their comments and
suggestions on an earlier draft of this paper.  The main results in
this paper were announced in a minicourse by the first author at a
conference at P.U.C.~in Rio de Janeiro in August 2001.  The first
author would like to thank the P.U.C. for their hospitality.

\section{Background on essential laminations}\label{essential_laminations}

In this section, we summarize the definitions and some of the basic
results about taut foliations and essential laminations (for more
detail see the references in \cite{Gabai1997prob}).  Throughout, $M$
will be a closed 3-manifold.  We'll begin with taut foliations:
\begin{defn}
  Let $\F$ be a foliation of $M$ by surfaces.  The foliation $\F$ is
  \emph{taut} if there is a loop in $M$ transverse to $\F$ which
  intersects every leaf.
\end{defn}
We will use $\til{\F}$ to denote the induced foliation of the
universal cover of $M$.  The \emph{leaf space} $L$ of $\til{\F}$ is a
simply-connected, not necessarily Hausdorff, $1$-manifold.  Some basic
properties of manifolds with taut foliations are summarized in the following
theorem:
\begin{Theorem}[Novikov, Rosenberg, Reeb]
  Let $M$ be a $3$-manifold with a taut foliation $\F$. Then either
  $M$ is finitely covered by $S^2 \times S^1$, and $\F$ is finitely
  covered by the product foliation by spheres, or $\til{M} = \R^3$
  foliated by planes. In particular, every leaf of $\F$ is
  incompressible, and every loop transverse to $\F$ is homotopically
  essential.
\end{Theorem}

A codimension one \emph{lamination} of $M$ is a foliation of a closed
subset, i.e.~a closed union of complete embedded surfaces---the {\em leaves} of
the lamination---which
intersect small open charts of $M$ in products of the form
$\R^2 \times K$, where $K$ is any closed subset of an interval, and each
leaf locally intersects the open set as a horizontal slice
$\R^2 \times \text{point}$. Globally, of course, a leaf might intersect
a product chart in infinitely many horizontal slices.
Typically, $K$ consists of a
union of intervals, isolated points, and Cantor sets.  Generalizing
the notion of tautness for foliations is:
\begin{defn}
  A lamination $\Lambda$ is \emph{essential} if each complementary
  region is irreducible and has boundary which is incompressible and
  end-incom\-pressible.  Moreover, we require that no leaf of $\Lambda$
  is a 2-sphere, or a torus leaf bounding a Reeb component.
\end{defn}
Here a leaf is end-compressible if it has an end-compressing
monogon---that is, a monogon $D$ properly embedded in a complementary
region $C$  which is not homotopic (rel boundary) into $\bdry C$.
An end-compressing monogon is an obstruction to finding a metric with
respect to which leaves of $\Lambda$ are minimal surfaces; excluding
such monogons ensures that essential laminations share many properties
in common with taut foliations.  In particular, manifolds with
essential laminations have universal cover $\R^3$ and any tight
transversal is homotopically essential (see \cite{GabaiOertel89}).

For an essential lamination, we can define the leaf space $L$ of
$\til{\Lambda}$, but this is now just an order tree, not a 1-manifold;
the vertices of $L$ are the set of non-boundary leaves of
$\til{\Lambda}$ together with the set of closed complementary regions.
If this leaf space is Hausdorff, we say $\Lambda$ is \emph{tight}.
This is equivalent to the leaves of $\til{\Lambda}$ being uniformly
properly embedded in $\til{M}$. That is, there is a proper increasing
function $f: \R^+ \to \R^+$ such that if $p$ and $q$ are distance $t$
apart in the path metric on $\til{\Lambda}$, then $p$ and $q$ are at
least distance $f(t)$ apart in $\til{M}$.  When $\Lambda$ is a
foliation, tightness is equivalent to the leaf space being $\R$.

Regardless of whether $\Lambda$ is tight, we have the following simple
fact:
\begin{Lemma}\label{uniform_constant}
  Let $\Lambda$ be an essential lamination. Then there is an
  $\epsilon>0$ such that every leaf $\lambda$ of $\til{\Lambda}$ is
  quasi-isometrically embedded in its $\epsilon$-neighbor\-hood.  Such
  an $\epsilon$ is called a \emph{separation constant} for $\Lambda$.
\end{Lemma}
\begin{proof}
  By compactness of $M$, we can pick an $\epsilon$ so that every ball
  of radius $2 \epsilon$ in $\til{\Lambda}$ is contained in a product
  chart.  Here, by ``product chart'' we also require that every
  transverse arc without backtracking and with
  endpoints in $\Lambda$ is tight, that is, can't
  be homotoped rel endpoints into a leaf.  Note that in $\til{M}$ any
  leaf $\lambda$ intersects a product chart at most once, as otherwise
  we can build a tight transversal loop to $\til{\Lambda}$ by starting
  with a transversal in the chart and an arc in $\lambda$.
  
  Now consider the $\epsilon$-neighborhood $N_\epsilon$ of a leaf
  $\lambda$ in $\til{\Lambda}$.  We can cover the larger neighborhood
  $N_{2 \epsilon}$ by a product charts .  As noted above, $\lambda$
  intersects each of these charts only once.  Because $M$ is compact,
  each of these charts has bounded geometry.  Therefore, $\lambda$ is
  quasi-isometrically embedded in $N_\epsilon$.  
\halmos \end{proof}

\subsection{Complementary regions}
Let $C$ be a complementary region to an essential lamination
$\Lambda$.  We will consider decompositions of $C$ along a union of
properly embedded annuli into two kinds of pieces, \emph{guts} and
\emph{interstitial regions}.  By definition, gut regions are compact
and interstitial regions are $I$-bundles over (typically noncompact)
surfaces.  The annuli of the decomposition are called the
\emph{interstitial annuli}, and we require that no two of them are
isotopic in a gut region.  By convention, if $C$ is itself an
$I$-bundle over a closed surface, the only allowed decomposition is
the trivial one with all of $C$ being the interstitial region.

Such a partition of $C$ into guts and interstices need not be unique:
if $R$ is an interstitial region $R = \Sigma \times I$ over a surface
$\Sigma$, and if $\Sigma_0 \subset \Sigma$ is a compact subsurface of
$\Sigma$ (possibly with boundary), then the subset $\Sigma_0 \times I$
can be included into the gut. Conversely, an $I$-bundle over a compact
surface which is contained in the gut can be included into the
interstitial region.  However, there exist several canonical (up to
isotopy) partitions of $C$ into guts and interstices. One such
canonical partition has the property that every interstitial region is
an $I$-bundle over a non-compact surface (see \cite{GabaiOertel89} or
\cite{GabaiKazez98} ).

Often, it's useful to make the distinction between essential
laminations and genuine ones:
\begin{defn}
  An essential lamination is {\em genuine} if some complementary
  region is not an $I$-bundle.  Equivalently,
  there is a decomposition of the complementary regions with non-empty
  guts.  
\end{defn}
An essential lamination which is not genuine can be filled in to a
foliation by adding leaves in the complementary $I$-bundles.

One way to think of the qualities essential and genuine for a
lamination is in terms of the proper essential surfaces (with
boundary) contained in a complementary region.  A lamination is
essential if the complement admits no essential sphere, disk or monogon;  that
is, it contains no essential surface of positive Euler characteristic.
An essential lamination is genuine if some complementary region
contains an essential surface of {\em negative} Euler
characteristic.

The following is an example of a genuine lamination which is a good
example to think about throughout this paper.
\begin{exa}\label{surface_bundle}
  Let $M$ be a surface bundle over the circle with fiber $\Sigma$ and
  pseudo-Anosov monodromy $\phi:\Sigma \to \Sigma$.  There are a pair
  of geodesic laminations of $\Sigma$ invariant under $\phi$, the
  stable and unstable laminations $\lambda^\pm$.  Since these are
  invariant, they suspend in the mapping torus $M$ to a pair of
  laminations $\Lambda^\pm$ of $M$ transverse to the surface
  fibration. The complementary regions of the geodesic laminations
  $\lambda^\pm$ of $\Sigma$ suspend to complementary regions of
  $\Lambda^\pm$. The complementary regions in $\Sigma$ are finite
  sided ideal polygons, each with at least $3$ sides. These suspend to
  ideal polygon bundles over $S^1$.  Such complementary regions
  decompose into compact neutered ideal polygon bundles over $S^1$
  (the guts) and cusps $S^1 \times \R^+ \times I$ (the interstitial
  regions). In particular, these laminations are genuine.  In
  $\til{M}$, the pair $(\til{M}, \til{\Lambda}^\pm)$ is a product
  $(\til{\Sigma}, \til{\lambda}^\pm)\times \R$.  As $\lambda^\pm$ is a
  geodesic lamination, the leaf space of $\til{\lambda}^\pm$ is
  Hausdorff.  Thus the leaf space of $\til{\Lambda}^\pm$ is also
  Hausdorff, and $\Lambda^\pm$ is tight.
\end{exa}

\subsection{Tangential geometry}

The tangential geometry of an essential lamination is controlled by
the following theorem of Candel \cite{Candel93}, which can be thought of as a
uniformization theorem for Riemann surface laminations:

\begin{Theorem}[Candel]\label{Candel_uniformize}
  Let $\Lambda$ be a Riemann surface lamination such that
  for every transverse measure $\mu$ we have:
  \[
  \chi(\mu)<0.
  \]
  Then there is a continuously varying leafwise metric on $\Lambda$
  where the leaves are locally isometric to $\H^2$.
\end{Theorem}
A lamination satisfying the conclusion of Candel's Theorem is said to
have \emph{hyperbolic leaves}.  For an essential lamination $\Lambda$
of an atoroidal $3$-manifold, every transverse measure has negative
Euler characteristic, and so $\Lambda$ has hyperbolic leaves.

Now we restrict attention to a taut foliation $\F$.  Candel's Theorem
is very useful in understanding the following bundle:
\begin{defn}
Let $\F$ be a taut foliation. The {\em circle bundle at infinity}, $E_\infty$,
is the circle bundle over $L$ whose fiber over a leaf $\lambda$ is the
circle $S^1_\infty(\lambda)$.
\end{defn}

While the bundle $E_\infty$ is easy to understand as a set, its
topology requires some discussion. We topologize it as follows.  For
each transversal $\tau$ to $\til{\F}$, the projection of $\tau$ to $L$
is an embedding, so it makes sense to talk about the restriction
$E_\infty|_\tau$.  Consider the restriction to $\tau$ of the unit
tangent bundle of $\F$, denoted $\mathrm{UT}\til{\F}|_\tau$.  For a
point $p \in \lambda$, every vector $v \in \mathrm{UT}_p\til{\F}$
determines a point $e(v) \in S^1_\infty(\lambda)$ by taking the
endpoint of the geodesic ray in $\lambda$ starting at $p$ in direction
$v$. Thus there is a bijection
\[
e_\tau \maps \mathrm{UT}\til{\F}|_\tau \to E_\infty|_\tau.
\]
We topologize $E_\infty$ by declaring that this map is a homeomorphism for each
$\tau$.

We verify that this is well-defined. Suppose $\tau_1$ and $\tau_2$ are
two transversals with the same projection to $L$. For each leaf
$\lambda$, the map
\[
e^{-1}_{\tau_2} \circ e_{\tau_1} \maps \mathrm{UT}_{\tau_1 \cap \lambda}\lambda \to \mathrm{UT}_{\tau_2 \cap \lambda}\lambda
\]
is determined by the geometry of a compact disk containing the points
$\tau_1 \cap \lambda$ and $\tau_2 \cap \lambda$. In particular, since
the geometry of the leaves $\lambda_t$ intersecting the $\tau_i$
varies continuously {\em on compact subsets}, the map
\[
e^{-1}_{\tau_2} \circ e_{\tau_1} \maps \mathrm{UT}\til{\F}|_{\tau_1} \to
\mathrm{UT}\til{\F}|_{\tau_2}
\]
is a homeomorphism.   Thus the topology of $E_\infty$ is well-defined.

This may seem like a tedious verification, but there actually is a
subtle point here. This construction uses in an essential way the full
power of Candel's Theorem. If we merely knew that the leaves of
$\til{\F}$ were uniformly coarsely quasi-isometric to $\H^2$, the map
$e_{\tau_2}^{-1} \circ e_{\tau_1}$ between unit tangent circles would
depend on the \emph{global} geometry of the leaf $\lambda$, and not
merely on the local geometry. But in general, the leaves $\lambda$ do
not vary continuously as subsets of $\til{M}$ in any reasonable sense.
This is the very subtlety that makes Candel's theorem nontrivial.
There is an alternate construction of the topology of $E_\infty$
without invoking Candel's theorem, but it is more involved, and uses
the Leaf Pocket Theorem~\ref{leaf_pocket_theorem}.

\section{Tight laminations}\label{tight_laminations}

In this section, we prove that certain kinds of tight genuine
laminations give rise to faithful actions on the circle.  Recall that
a genuine lamination is very full if the complementary regions are all
finite-sided ideal polygon bundles over $S^1$.  Our basic result here
is:
\begin{Theorem}\label{order_tight_full}
  Let $M$ be an orientable atoroidal $3$-manifold containing a very
  full essential lamination $\Lambda$ which is tight. Then there is a
  faithful representation
  \[
  \rho:\pi_1(M) \to \homeo^+(S^1).
  \]
  Moreover, the image of $\rho$ preserves a dense lamination of $S^1$.
\end{Theorem}
Here, a lamination $\lambda$ of $S^1$ is essentially just a geodesic
lamination of $\H^2$.  Intrinsically, $\lambda$ is a closed, symmetric subset of
$\lbrace S^1 \times S^1 - \text{diagonal} \rbrace$ so that no two
pairs of points $(p_1, p_2)$ and $(q_1, q_2)$ in $\lambda$ are linked on $S^1$.

The hypothesis of very full is much stronger than just having solid
torus guts; the former requires in addition that the interstitial
regions have as little topology as possible.  An typical example of a
lamination with solid torus guts which is not very full is depicted in
Figure~\ref{splitting_C}.   However, the Filling
Lemma~\ref{filling_lemma} shows that a lamination with solid torus
guts can be included as a sublamination of a very full lamination.
Moreover, this process preserves tightness.  So as an immediate
consequence of the Filling Lemma and Theorem~\ref{order_tight_full} we
have:
\begin{Theorem}\label{main_thm}
Let $M$ be an orientable atoroidal $3$-manifold, and $\Lambda$ a tight lamination
with solid torus guts. Then $\pi_1(M)$ has a faithful representation into
$\homeo^+(S^1)$.
\end{Theorem}

While the restriction of solid torus guts might initially seem quite
strong, in fact many constructions of essential laminations yield
examples of this type.  Moreover, for non-Haken manifolds, the gut
regions are always handlebodies of some sort \cite{Brittenham93}.

Now let's prove Theorem~\ref{order_tight_full}.

\begin{proof}[Proof of Theorem~\ref{order_tight_full}]
  For any essential lamination $\Lambda$, the universal cover
  $\til{M}$ is topologically $\R^3$, and the lifted lamination
  $\til{\Lambda}$ is topologically a product of a lamination $\lambda$
  of $\R^2$ with $\R$ \cite{GabaiKazez1997GA}.  Moreover, the leaves
  of $\lambda$ are proper.  Thus $\lambda$ is quite close to being a
  geodesic lamination of $\H^2$.   As our $\Lambda$ is very full, the
  complementary regions of $\lambda$ are all finite-sided ideal
  polygons.  In essence, we'll construct an action of $\pi_1(M)$ on
  $(\R^2, \lambda)$ so that the map $(\R^3, \til{\Lambda}) \to (\R^2,
  \lambda)$ is equivariant, and from this we'll get the required
  action on a circle.  As $\Lambda$ is very full, the complementary
  regions of $\til{\Lambda}$ are all products $P_i \times \R$ for
  finite sided ideal polygons $P_i$.  When turning $(\R^3,
  \til{\Lambda})$ into a product, there are essentially two possible
  ways of flattening each $P_i \times \R$ into a complementary region
  of $(\R^2, \lambda)$.  The key to constructing the desired action on
  $S^1$ is to do this flattening in a consistent, equivariant way.
  
  So let's begin the actual construction, which involves fitting
  the $P_i$ together as ideal polygons in the unit disk.  First,
  downstairs in $M$, choose an orientation on the core curve of each
  complementary region of $\Lambda$.  Lifting upstairs, this defines a
  (dual) orientation on each $P_i$ which is preserved by the action of
  $\pi_1(M)$.  As with any essential lamination, there is a natural
  leaf space $L$ associated to $\til{\Lambda}$.  The leaf space $L$ is
  an order tree associated to $\til{\Lambda}$ whose vertices are the 
  set of non-boundary leaves of $\til{\Lambda}$ together with the set
  of closed complementary regions. This order tree can be canonically included
  in an $\R$-order tree (see e.g.~\cite{GabaiKazez1997MRL}); this makes the
  explanation of the construction easier to follow, so we suppose this
  has been done. As $\Lambda$ is tight, $L$ is Hausdorff.
  
  We will insert the $P_i$ as ideal polygons in the unit disk $D$ via
  a map $\pi:P_i \to D$ which satisfies the following conditions:
  \begin{enumerate}
  \item If $i \ne j$, then the closed polygons $\pi(P_i)$ and
      $\pi(P_j)$ are disjoint.
  \item The circular ordering on the vertices of $\pi(P_i)$
      induced by $\partial D$ agrees with that arising from
      the orientation on $P_i$.
  \item Suppose  $P_i$ is joined to $P_j$ by an embedded arc in $L$ joining
      a side $e_i$ of $P_i$ to a side $e_j$ of $P_j$.  Then the corresponding edges
      $\pi(e_i)$ and $\pi(e_j)$  separate the interiors of $\pi(P_i)$ and $\pi(P_j)$
      in $D$.
  \end{enumerate}
  Here's how to construct $\pi$.  Because the $\R$-order tree $L$ is
  Hausdorff, it is the union of a countable number of finite
  simplicial trees $L_j$ \cite[\S 3]{GabaiKazez1997MRL}.  Here, $L_j$
  is subdivided and extended to get $L_{j+1}$.  Consider $L_j$ as a
  simplicial complex with no edges of valence two.  Each interior
  vertex of $L_j$ corresponds to some complementary region of
  $\til{\Lambda}$, and thus to some $P_i$.  The edges of $L_j$ coming
  into a vertex correspond to the edges of the associated $P_i$.
  Thus, the orientation on $P_i$ induces a cyclic ordering of the
  edges around each vertex of $L_j$.  There is a unique embedding
  $\pi$ of $L_j$ into the interior of $D$ which respects these
  cyclic orderings (i.e. at each vertex $v$ of $L_j$ the ordering of
  edges agrees with the clockwise ordering of the image edges about
  $\pi(v)$).  For each vertex $v$ in $L_j$, place the corresponding
  $P_i$ as an ideal polygon in $D$ so that it contains $\pi(v)$.  The
  resulting polygons satisfy (1-3), and moreover their placement
  was unique up to homeomorphism of $D$.  Inducting up the $L_j$ gives
  us the required map $\pi$.

  Now, the closure 
  \[
  K = \overline{\bigcup_i\pi(P_i)}
  \]
  might not be all of $D$, but the components of $D-K$ are easy to
  understand.  Any component $C$ of $D-K$ is a polygon whose sides
  alternate between arcs of $\bdry D$ and geodesics in $D$ which are
  limits of edges in $\pi(\cup_i \bdry P_i)$.  Moreover, we claim $C$ must be a
  quadrilateral.  To see this, pick one $P_i$ in each component of $D
  - C$ and join them by a finite simplicial tree $L'$ in $L$.
  Embedding $L'$ in $D$, we see that if there were more than two
  components to $D-C$ , then $L'$ would have to have a vertex in $C$.
  
  If we quotient out $\partial D$ by collapsing the arcs in $\partial
  D - K$ to points, we get a circle because $\til{\Lambda}$ has no
  isolated leaves.  Let $\phi:\partial D \to S^1$ be this monotone
  quotient map.  The closure of the edges of the $\pi(\partial P_i)$
  is a geodesic lamination $\lambda_K$ of $D$.  The map $\phi$ gives a
  map on pairs of endpoints of geodesics of $\lambda_K$ which
  preserves unlinking, so there is a well-defined image lamination
  $\lambda = \phi_*(\lambda_K)$ in $S^1$.  The endpoints of leaves of
  $\lambda$ are dense in $S^1$.  The action of $\pi_1(M)$
  on $\bigcup_i P_i$ induces an action on $\lambda$.  Going back to
  our original construction of $\pi$, it is not hard to see that this
  extends to a continuous action on $S^1$ which preserves
  orientation.
  
  To finish the proof, we need show this action is faithful.  Consider
  a nontrivial element $\alpha \in \pi_1(M)$ which acts trivially on
  $S^1$.  Then the action of $\alpha$ on $\til{M}$ would stabilize
  every complementary region.  So $\alpha$ defines the core curve of
  every complementary region.  It follows that $\alpha$ is central in
  $\pi_1(M)$, and that $\pi_1(M)$ contains a $\Z + \Z$.  This violates our
  assumption that $M$ is atoroidal.  Thus we have constructed the
  required faithful action on a circle.
\halmos \end{proof}

\subsection{Flips}\label{flips_one}
If there are $n$ complementary regions to $\Lambda$, the above
construction gives $2^n$ representations which are pairwise
non-conjugate in $\homeo^+(S^1)$ (they fall into only $2^{n-1}$
conjugacy classes in the full group $\homeo(S^1)$).
  
  Abstractly, this reflects an interesting operation on
  representations of groups
  \[
  \sigma:G \to \homeo^+(S^1)
  \]
  whose image preserves a lamination $\lambda$ of $S^1$.  Given an orbit class of
  complementary regions to $\lambda$, we get a new representation as
  follows.  Let $P$ be a complementary region of $\lambda$.  We can
  ``cut'' $D$ along the edges of each translate $\sigma(\gamma)(P)$
  into pieces which are $\sigma(\gamma)(P)$ itself, and the pieces on
  the other sides of the boundary edges of $\sigma(\gamma)(P)$.  Then
  we can reverse the orientation of $\sigma(\gamma)(P)$ and glue the
  other pieces back so that each pair of sides is glued in the same
  way as before, but the relative orientations are now switched.  Do
  this in some order for each $\gamma \in G$, to get a new lamination
  $\lambda^P$ on which $G$ acts.  This action extends in the obvious
  way to a representation $\sigma^P:G \to \homeo^+(S^1)$.  Note that
  $(\sigma^P)^P = \sigma$ and $(\lambda^P)^P = \lambda$ for any $P$.
  We call this operation {\em the flip of $\sigma$ along $P$}. 
  
  If $\lambda$ only has one orbit class of complementary region, this
  operation merely reverses the orientation of $S^1$, but in general
  it seems somewhat mysterious.  It might be interesting to analyze
  this operation for some familiar $G$. For instance, one could look
  at $G = \pi_1(M)$ for $M$ a hyperbolic $3$-manifold fibering over
  $\Sigma$, and $\lambda$ the lift of the stable lamination of the
  monodromy of the fibering.

  The moral of these examples is that if a group $G$ acts on a ($\R$-,
  order-) tree $T$ and there is at least one way of embedding $T$ in
  the plane in a $G$-equivariant way, there are usually many such
  ways.
  
  In the sequel, we will be interested in obtaining representations of
  $G = \pi_1(M)$ in $\homeo(S^1)$ with Euler class $e \in H^2(G)$
  equal to $0$ (see Section~\ref{sec_orderings}).  We suspect the flip
  operation might be a useful tool in this regard. It is not hard to
  see that if $\lambda$ has $n$ equivalence classes of complementary
  polygons, there is an expression of $e$ as a sum
  \[
  e = \sum_{i=1}^n e_i
  \]
  such that the Euler class of $\sigma^{P_j}$ is 
  \[
  e_{P_j} = \left( \sum_{i=1}^n e_i \right) - 2e_j.
  \]

\subsection{Laminations and cataclysms}\label{sec_non_tight}

The rest of this section is devoted to generalizations of
Theorem~\ref{main_thm} where we partially relax the requirement that
the lamination be tight.  To start, we need to get a handle on what a
non-tight lamination looks like. 

If $\Lambda$ is not tight, by definition the leaf space $L$ of
$\til{\Lambda}$ is a non-Hausdorff order tree.  Recall that two leaves
$\mu$ and $\lambda$ in $\til{\Lambda}$ are comparable if they can be
jointed by a tight transversal, and incomparable otherwise.  A
sequence $\{\mu_i\}$ of leaves in $L$ is \emph{monotone ordered} if
all of the $\{\mu_i\}$ lie in an ordered segment $I \subset L$, so
that the $\mu_i$ form an increasing sequence in $I$.  We define a
\emph{cataclysm} to be a collection of incomparable leaves $\lbrace
\lambda_j \rbrace$ of $\til{\Lambda}$, called the \emph{limit leaves},
for which there is a monotone ordered sequence $\mu_i$,
called the \emph{approximating sequence}, which converges on compact
subsets in the Hausdorff topology to the union of the $\lambda_j$;
that is, for all compact subsets $K \subset \til{M}$,
\[
\lim_{i \to \infty} \left( \mu_i \cap K \right) = \left( \bigcup_j  \lambda_j \right) \cap K
\]
in the Hausdorff topology on closed subsets of $K$.

For any sequence $p_i \in \mu_i$ which converges to some $p \in
\lambda_j$, the sequence of pointed metric spaces $(\mu_i,p)$
converges to $(\lambda_j,p)$. In particular, if $p_i$ and $q_i$ are
two sequences with $p_i,q_i \in \mu_i$ where $p_i \to p \in \lambda_j$
and $q_i \to q \in \lambda_k$ for $\lambda_j \ne \lambda_k$, then the
leafwise distance from $p_i$ to $q_i$ in $\mu_i$ goes to infinity.

We define an equivalence relation on cataclysms as follows: two
cataclysms $C = (\lbrace \lambda_j \rbrace,\lbrace \mu_i \rbrace)$ and
$C' = (\lbrace \lambda_{j'}' \rbrace, \lbrace \mu_{i'}' \rbrace)$ are {\em
  equivalent} if $\lbrace \lambda_j \rbrace$ and $\lbrace \lambda_{j'}'
\rbrace$ are equal as {\em sets}, and if for some $N$, the union
\[
\lbrace \mu_j \rbrace_{j\ge N} \cup \lbrace \mu_{j'}' \rbrace_{j'\ge N}
\]
can be totally ordered in $L$.

The equivalence class of a cataclysm is determined by the collection
$\lbrace \lambda_j \rbrace$ of leaves, since for any $p \in \lambda_j$
and any $p_i \to p$ contained in $\til{\Lambda}$ on the side of
$\lambda_i$ which contains the other $\lambda_j$, the leaves $\mu_i$
with $p_i \in \mu_i$ are an approximating sequence, unique up to
equivalence. It follows that the stabilizer of a cataclysm acts by
permutations on the set of limit leaves, and by isometry on their
union, thought of as a subset of $\til{M}$.

The following condition is natural for trying to build circle actions:
\begin{defn}
  A lamination $\Lambda$ has \emph{orderable cataclysms} if for each
  equivalence class of cataclysm $[C]$, there is an ordering on the set of
  limit leaves $\lbrace \lambda_i \rbrace$ of $[C]$ which is invariant
  under the action of the stabilizer of $[C]$ in $\pi_1(M)$.
\end{defn}
We will show that orderable cataclysms can replace tight as a
hypothesis in Theorem~\ref{main_thm}, but let's first give an example
of this condition.

\begin{exa}\label{orderable_anosov}
  Let $X$ be a pseudo-Anosov flow, and denote the two dimensional
  stable and unstable singular foliations $\F^s$ and $\F^u$.  For ease
  of notation, we will focus on $\F^u$ but of course similar
  statements are true for $\F^s$.  The foliation $\F^u$ lifts to a
  singular foliation $\til{\F}^u$ of $\til{M}$ invariant under the
  lifted flow $\til{X}$. Each nonsingular leaf of $\til{\F}^u$ is
  topologically a plane foliated by flow lines of $\til{X}$. If
  $\lambda$ is a nonsingular leaf of $\til{\F}^u$, the flow along
  $\lambda$ by $\til{X}$ compresses the flowlines together. Thus the
  holonomy along this foliation is defined in the positive direction
  for all time. In particular, it is easy to see from this that
  $\lambda$ is isometric to $\H^2$.  The foliation $\lambda \cap
  \til{X}$ is a foliation of $\H^2$ by geodesics asymptotic to a
  unique point in $S^1_\infty$. A singular leaf $\lambda$ of
  $\til{\F}^u$ has the following structure.  Let $h_1,\dots,h_n$ be
  copies of $\H^2$ foliated by geodesics asymptotic to $p_i \in
  S^1_\infty(h_i)$, and let $\gamma_i \subset h_i$ be one of these
  geodesics.  Then $\gamma_i$ separates $h_i$ into two sides,
  $h_i^\pm$. We obtain a quotient space
  \[
  l = \bigcup_i h_i
  \]
  by gluing $h_i^+$ to $h_{i+1}^-$ via an isometry taking $\gamma_i$
  to $\gamma_{i+1}$, where the indices $i$ are taken modulo $n$. Then
  $\lambda = l$, and the $h_i$ are called the \emph{faces} of the
  singular leaf $\lambda$. The image of the $\gamma_i$ in $l$ is a
  singular orbit of $\til{X}$ and covers a circular singular orbit in
  $M$.
  
  The singular foliation $\F^u$ can be split open into a lamination
  $\Lambda^{u}$ whose leaves are exactly the nonsingular leaves of the
  $\til{\F}^{u}$ together with one leaf for every face of a singular
  leaf.  If $X$ is pseudo-Anosov but not Anosov, then $\Lambda^{u}$ is
  a very full genuine lamination, and if $X$ is Anosov then
  $\Lambda^{u} = \F^{u}$.  The foliation of the leaves of $\F^{u}$ by
  flow lines of $X$ gives rise to foliations of the leaves of
  $\Lambda^{u}$.  Consider a cataclysm $C$ with limit leaves
  $\lambda_i$ and approximating sequence $\mu_j$.  The foliations of
  the $\mu_j$ by flowlines of $\til{X}$ converge on compact sets to
  the foliations of the $\lambda_i$.  So if we pick points $p_{ij} \in
  \mu_j$ so that $p_{ij}$ is very close to $\lambda_i$, there is a
  natural order on the $p_{ij}$ coming from the order on the leaf space
  of the $\til{X}$ foliation of $\mu_j$, which is $\R$. This determines
  an order on the indices $i$.  As $j \to \infty$, the order on the
  indices $i$ is eventually constant, and therefore determines a
  natural order on the set $\lbrace \lambda_i \rbrace$, which is
  invariant under the stabilizer in $\pi_1(M)$ of $[C]$.  Thus
  $\Lambda^u$ has orderable cataclysms.
  
  This natural order structure on the cataclysms of pseudo-Anosov
  flows was observed by Fenley in \cite{Fenley98}, where he gives
  examples of pseudo-Anosov flows where $\Lambda^u$ is not tight.
\end{exa}

\begin{exa}
  Not all cataclysms are orderable.  Here is a simple example.  Start
  with a tight lamination $\Lambda$ with a complementary region which
  is an ideal square bundle over the circle where the monodromy is the
  rotation through angle $\pi$.  Then fill in that region with a
  saddle bundle over $S^1$. In the universal cover, the bundle of
  saddles lifts to $\R$ leaves which limit in either direction to a
  pair of distinct leaves corresponding to opposite sides of the ideal
  square.   Either end gives a cataclysm with two limit leaves.  The
  stabilizer of the cataclysm corresponds to the monodromy of the core
  of the complementary region downstairs; in particular, it
  interchanges these two limit leaves, and this cataclysm is not
  orderable.
\end{exa}

Now we'll show:
\begin{Theorem}\label{thm_ord_cataclysms}
  Let $M$ be an orientable atoroidal $3$-manifold containing a
  lamination $\Lambda$ with solid torus guts and orderable cataclysms.
  Then $\pi_1(M)$ has a faithful representation in $\homeo^+(S^1)$.
\end{Theorem}
\begin{proof}
  By the Filling Lemma~\ref{filling_lemma}, $\Lambda$ can be filled to
  a very full lamination.  The proof that filling preserves tightness
  further shows that filling preserves the set of equivalence classes
  of cataclysms, and that no added leaf is a limit leaf of a
  cataclysm. In particular, the filled lamination has orderable
  cataclysms.  So from now on, we'll assume that $\Lambda$ is very
  full.
  
  The proof is now morally the same as for
  Theorem~\ref{order_tight_full}, but this time it's easier to do
  things a little more abstractly.  A circular order on a set $S$ is
  an assignment to each triple of distinct elements $(s_1, s_2, s_3)$
  an orientation of clockwise or anti-clockwise, satisfying the
  obvious rules coming from triples of points on $S^1$.  Here, we will
  show that the set of ends $E$ of $L$ has a natural circular order
  which is invariant under the action of $\pi_1(M)$.  We will complete
  $E$ with respect to this circular order to get the needed circle.
  (For a detailed definitions and basic properties of circular orders
  on sets and left-invariant circular orders on sets acted on by a
  group, see e.g.~\cite{Tararin2001}.)

  By definition, a circular order on an order tree
  $L$ consists of two things: a circular order on sets of segments
  with only a vertex in common, and a linear order on sets of segments
  which differ only by a vertex (for details see
  \cite{GabaiKazez1997MRL}).  For the leaf space $L$ of our lamination
  $\til{\Lambda},$ we can define a circular order which is invariant
  under the action of $\pi_1(M)$ as follows.  A set of segments with
  only a vertex in common corresponds to the set of faces of a closed
  complementary region $C = P_i \times \R$ of $\til{\Lambda}$.  As
  before, if we fix orientations of the core curves of the
  complementary regions $M - \Lambda$, the faces of $C$ have a natural
  circular order.  A set of segments which differ only in a vertex
  corresponds to a cataclysm, and so has an order by assumption.  In
  the cataclysm case, we can certainly choose our orderings to be
  $\pi_1(M)$-equivariant.
  
  Thus $L$ is an order tree with a $\pi_1(M)$-invariant circular
  order.  Let $E$ be the set of ends of $L$.  If $(e_1, e_2, e_3)$ are
  three distinct ends in $E$, we define their circular order as
  follows.  Take any point $p$ in $L$ and consider rays $r_i$ starting
  at $p$ and ending at $e_i$.  By looking at the (circularly ordered)
  subtree $T = \cup r_i$, we get a circular order of the ends of the
  $r_i$.  This order is independent of the choice of point $p$, and so
  we have a natural circular order on $E$ which is invariant under the
  action of $\pi_1(M)$.  Up to homeomorphism, there is a unique way to
  embed $E$ into a circle $S^1$ so that the embedding respects the
  circular orderings and is continuous with respect to the topology on
  $E$ induced by the order.  The closure of the image of $E$ may
  omit some gaps, which we can collapse to get the promised circle on
  which $\pi_1(M)$ acts.  As in the proof of
  Theorem~\ref{order_tight_full}, if some element $\alpha \in
  \pi_1(M)$ acts trivially on this circle, it must fix every
  complementary domain; then $\pi_1(M)$ has a nontrivial center, violating the
  assumption that $M$ is atoroidal.  
\halmos \end{proof}

\begin{Corollary}\label{cor_anosov_flow}
  Let $M$ be an atoroidal $3$-manifold which admits a pseudo-Anosov flow
  $X$. Then $\pi_1(M)$ admits a faithful representation in
  $\homeo(S^1)$.
\end{Corollary}
\begin{proof}
  If $X$ is pseudo-Anosov but not Anosov, the singular unstable
  foliation can be split open to a very full genuine lamination. As
  explained in Example~\ref{orderable_anosov}, these laminations have
  orderable cataclysms and so Theorem~\ref{thm_ord_cataclysms}
  applies.  If $X$ is Anosov, the unstable foliation is a taut
  foliation, and the circle action comes from
  Theorem~\ref{universal_circle_action_faithful}.
\halmos \end{proof}

\subsection{Invariant laminations for cataclysms}

For the sake of completeness, we give a more precise
description of the non-Hausdorff behavior of the leaf
space at a cataclysm. The following theorem shows that the set of limit leaves
of a cataclysm can be parameterized by a geodesic lamination of $\H^2$ in a (topologically)
canonical way.

\begin{Theorem}\label{cataclysm_lamination}
  Let $\Lambda$ be an essential lamination, and let $C$ be a cataclysm
  of $\Lambda$. Then we can associate a lamination $L$ of $\H^2$ to
  $C$ in such a way that the complementary regions to $L$ are in
  $1$-$1$ correspondence with the limit leaves of $C$, and the
  stabilizer $\Gamma \subset \pi_1(M)$ of $[C]$ acts by homeomorphisms
  of $\H^2$ which preserve $L$, and permute the complementary regions
  by the permutation action on the limit leaves of $C$.
\end{Theorem}
\begin{proof}
  For each limit leaf $\lambda_i$, pick some point $p_i$, and for each
  approximating leaf $\mu_j$ pick $p_{ij}$ such that the $p_{ij}$
  limit to $p_i$ as $j \to \infty$, for each $i$.  Pick some small
  $\epsilon$, and consider for each limit leaf $\lambda_i$ and each
  approximating leaf $\mu_j$ the subspace $\mu_j(i) \subset \mu_j$
  consisting of points which are within $\epsilon$ of $\lambda_i$.
  Then for sufficiently small $\epsilon$, the sets $\mu_j(i)$ and
  $\mu_j(k)$ are disjoint for $i \ne k$, since for any lamination with
  bounded geometry, leaves are uniformly properly embedded in their
  $\epsilon$ neighborhoods for some $\epsilon$, and therefore pairs of
  points on incomparable leaves of $\til{\Lambda}$ are a uniform
  distance apart. On the other hand, by the Leaf Pocket
  Theorem~\ref{leaf_pocket_theorem}, for large $j$, $\mu_j(i)$
  contains a $\delta$-net of points in the circle at infinity
  $S^1_\infty(\mu_j)$ in the visual metric as seen from $p_{ij}$.  In
  particular, the boundary of the convex hull of the limit set of
  $\mu_j(i)$ defines a collection of geodesics which separates most of
  $\mu_j(i)$ from $\mu_j(k)$ for $i \ne k$ in some definite order; the
  closure of the union of such geodesics is a geodesic lamination
  $L_j$ of $\mu_j$. This lamination is dual to a planar order tree
  $T_j$ whose vertices are the $i$ for which $p_{ij}$ is sufficiently
  close to $\lambda_i$.  In particular, there are a sequence of
  inclusions $T_j \to T_{j+1} \to \dots$ and the union $T_\infty$ is a
  Hausdorff planar order tree dual to a lamination $L$ of $\H^2$,
  which is a limit of the $L_j$ under an appropriate sequence of
  homeomorphisms (\emph{not} isometries) from $\mu_j \to \H^2$. It is
  clear that $L$ does not depend on the choice of approximating
  leaves, since if we insert some $\mu_{j+1/2}$ between $\mu_j$
  and $\mu_{j+1}$ then there are inclusions $T_j \to T_{j+1/2}
  \to T_{j+1} \to \dots$ and the limit is unchanged.  In particular,
  $L$ depends only on $[C]$, and therefore the stabilizer of $[C]$
  acts on it by automorphisms.
\halmos \end{proof}

\section{The Filling Lemma}\label{section_filling_lemma}

This section is devoted to the proof of the following lemma, which is
of independent interest.  It resolves the disparity between having
solid torus guts and solid torus complementary regions, and is used in
Section~\ref{tight_laminations} to reduce the construction of actions
on circles to the latter case.

\begin{filling_lemma}\label{filling_lemma}
  Let $M$ be an orientable $3$-manifold and $\Lambda$ a genuine
  lamination such that for some decomposition of $M - \Lambda$ into
  interstices and guts, the guts are neutered ideal polygon bundles
  over $S^1$. Then $\Lambda$ can be filled to a genuine lamination
  $\overline{\Lambda} \supset \Lambda$, whose complementary regions
  are all ideal polygon bundles over $S^1$. That is,
  $\overline{\Lambda}$ is very full. Moreover, if $\Lambda$ is tight,
  so is $\overline{\Lambda}$.
\end{filling_lemma}

Before giving the proof, let us point out another application.  Gabai
and Kazez have shown that if an atoroidal manifold has a genuine
lamination with some complementary region a solid torus, then
homotopic homeomorphisms are isotopic \cite{GabaiKazez1997GA}.  The
Filling Lemma extends that result to manifolds with genuine
laminations where there is a decomposition of $M - \Lambda$ with a
solid torus gut region.   Now we'll prove the Filling Lemma.

\begin{proof}
  First, fill all product complementary regions with foliation.  Now,
  each complementary region $C$ of $\Lambda$ is a union of gut pieces
  $G_i$ together with a finite collection of interval bundles $J_k$
  (the bases of the $J_k$ need not be compact).  The $J_k$ are glued
  to the $G_i$ along the interstitial annuli.  The first step of the
  proof is to add leaves to reduce to the case where each $J_k$ is a
  trivial $I$-bundle over a surface with a single boundary component.
  Then, each $J_k$ is attached to a single $G_i$ along one interstitial
  annulus, and each $C$ contains exactly one gut region.  After this
  basic reduction, we'll add on a product foliation to the boundary
  leaves of $C$, and perturb it so that the new $C$ is an ideal polygon
  bundle over $S^1$.

  \subsection{Basic reduction}
  First, we'll insert finitely many leaves so that each
  $I$-bundle $J_k$ has only one interstitial annulus in its
  boundary.  Let $J_k$ be an $I$-bundle over a base with at least
  two boundary components.  Pick an interstitial annulus $A$ bounding
  $J_k$, joining it to a gut region $G_i$.  Let $\lambda$ and $\mu$
  be the leaves of $\Lambda$ that $A$ runs between.  
  \begin{figure}[h]
    \begin{center}
      \includegraphics[scale=0.5]{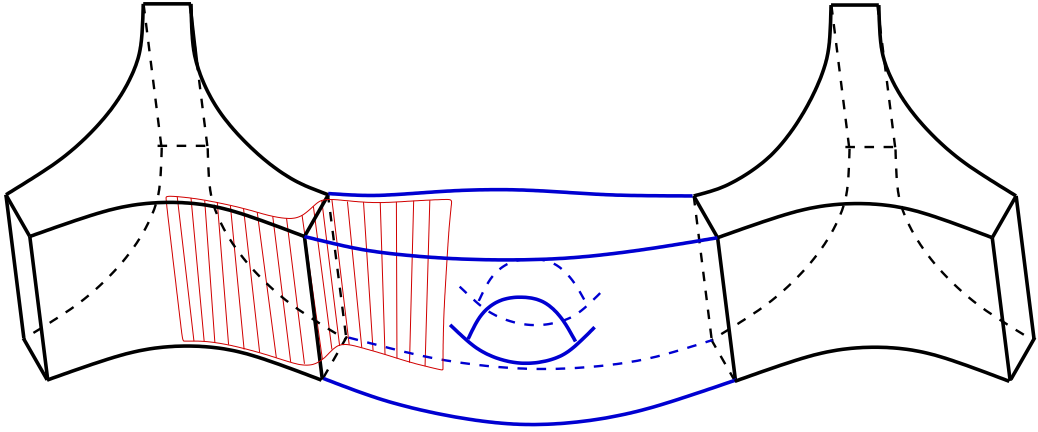}
    \end{center}
    \caption{      
      A $J_k$ with two boundary components joins two gut regions $G_i$
      and $G_j$.  The interstitial annulus $A$ we're considering is
      one on the left.  Also shown is part of the new leaf $\nu$ which
      we're adding.  The leaf $\nu$ is carried by the branched surface $B$.  } \label{splitting_C}
  \end{figure}
  
  Build a branched surface $B$ from $A$, $\mu$, and $\lambda$, bending
  $A$ as shown in Figure~\ref{splitting_C}.  The leaf we add will be
  carried by $B$, and will consist of one copy of $A$ glued to stuff
  parallel to $\mu$ and $\lambda$.  Explicitly, take countably many
  copies of the two leaves $\lambda$ and $\mu$ in a product
  neighborhood, which accumulate only along $\lambda$ and $\mu$.  Take
  the union of this with $A$ to get a two-complex which is singular
  along countably many circles accumulating to the two boundary
  circles of $A$.  Construct a surface carried by $B$ by performing a
  normal sum operation on each singular circle using a consistent
  orientation.  Let $\nu$ be the component of this new surface which
  contains most of $A$.  Add $\nu$ to $\Lambda$ to get $\Lambda'$.  To
  decompose the complement of $\Lambda'$ into $I$-bundles and guts,
  use essentially the same interstitial annuli as for $\Lambda$.  In
  the complement of $\Lambda'$, the $I$-bundle $J_k$ has been split
  into a $J'_k$ with one fewer boundary components, and another
  $I$-bundle with one boundary component (consult
  Figure~\ref{splitting_C}).  Doing this process at most once for each
  of our original interstitial annuli gives us a new $\Lambda$ where
  each $J_k$ has only one boundary component.
  
  Now, we'll add finitely many more leaves so that each $J_k$ is a
  trivial $I$-bundle.  Suppose some $I$-bundle $J_1$ is nontrivial,
  and let $G$ be the gut region it's attached to via the interstitial
  annulus $A_1$.  Let $\Sigma_1$ be a copy of the base of $J_1$
  sitting in the middle of $J_1$.  Let $J_2$ be an $I$-bundle attached
  to $G$ along $A_2$.  Choose a horizontal surface $\Sigma_2$ in $J_2$
  which has exactly one boundary component on $A_2$.  Let $\nu$ be the
  surface created by gluing $\Sigma_1$ and $\Sigma_2$ together via an
  annulus in $G$.  If we pick $J_2$ so that $A_1$ and $A_2$ are
  adjacent on $\bdry G$, the new lamination $\Lambda' = \Lambda \cup
  \mu$ has essentially the same guts as $\Lambda$.  If $J_2$ was also
  twisted, adding $\nu$ creates a new $J'_k$ with two boundary
  components, so in this case we add a second leaf as in the previous
  step so that all $I$-bundles still have a single boundary component.
  In any event, there are now fewer twisted $I$-bundles $J'_k$ in the
  complement of $\Lambda'$.  Thus, after adding finitely many leaves,
  we can make all of the $J_k$ trivial $I$-bundles over surfaces with
  one boundary component.

  \subsection{Main argument}   
  Consider one of our complementary regions $C$.  We need to fill
  $\Lambda$ so that $C$ becomes an ideal polygon bundle over $S^1$.
  As a rough outline, we'll begin by foliating each interstitial
  region $J_k$ by planes and annuli in a very standard way.  Here, by
  ``plane'' we really mean a disc with a part of its boundary removed,
  and similarly for annuli.  The induced foliation of each
  interstitial annulus will consist of a countable collection of
  circles, and a bunch of lines spiraling out to the circle leaves.
  Each annulus leaf has a single circle boundary component and many
  line boundary components.  We'll split open the foliation of each
  $J_k$ along one of the annulus leaves and glue these foliations
  together across the gut $G$ to enlarge our original lamination
  $\Lambda$.  If we glue with care, the new leaves are all either
  annuli or planes.  The new interstitial $I$-bundles will have bases
  which are annuli of the form $S^1 \times [0,1)$, and so the
  complementary regions are all ideal polygon bundle over $S^1$.
  
  We'll begin by constructing the foliation of each interstitial
  region.  Let $J$ be an interstitial region $\Sigma \times I$, where $\Sigma$
  is has exactly one boundary loop $\gamma$.  Let $A=\gamma \times I$.
  We will foliate $J$ by choosing some representation
  \[
  \til{\sigma} \maps \pi_1(\Sigma) \to \homeo(I)
  \]
  so that the induced foliation of $J$ as a foliated bundle has the
  above properties.  We will construct $\til{\sigma}$ from a
  representation $\sigma \maps \pi_1(\Sigma) \to \PSL{2}{\R}$ by
  lifting $\sigma$ to $\til{\PSL{2}{\R}}$.  The group
  $\til{\PSL{2}{\R}} \subset \til{\homeo}(S^1)$ acts on $\R$, and we
  can identify $\homeo(\R)$ with $\homeo(I)$.  We'll want the
  representation $\sigma$ to be faithful, and also that the
  nontrivial images are hyperbolic elements.  Moreover, we'll need
  that in the lift
  \[
  \til{\sigma}\maps \pi_1(\Sigma) \to \til{\PSL{2}{\R}}
  \]
  the rotation number of $\gamma$ is zero.  This way,
  $\til{\sigma}(\gamma)$ will have countably many fixed points, and
  will alternate between being increasing and decreasing on the
  complementary intervals.
  
  Here's why such a $\sigma$ exists.  The group $\pi_1(\Sigma)$ is a
  free group, possibly on infinitely many generators.  First, consider
  the case where $\pi_1(\Sigma)$ is finitely generated.  Note that the set of
  faithful representations is dense in the space of all $\PSL{2}{\R}$
  representations; this is because it is the complement of 
  \[
  \bigcup_{g \in \pi_1(\Sigma)} \setdef[\big]{\rho \maps \pi_1(\Sigma) \to
    \PSL{2}{\R}}{\rho(g) = I}
   \]
   which is a countable union of proper algebraic subvarieties.  For
   the trivial representation $\tau$, if we look at the trivial lift
   $\til{\tau}$ then the rotation number of $\til{\tau}(\gamma)$ is
   $0$.  For $\sigma$ near $\tau$ for which $\sigma(\gamma)$ is
   hyperbolic, the rotation number of $\til{\sigma}(\gamma)$ near
   $\til{\tau}(\gamma)$ is an integer arbitrarily close to $0$.  Thus
   it is $0$, and we can construct the needed $\sigma$.  Moreover, by
   taking a generic representation $\til{\sigma}$, we can also require
   that the stabilizer of all but a countable set of points in $\R$
   will be trivial, and the stabilizer of each element of that
   countable set will be isomorphic to $\Z$.
   
   Now suppose $\pi_1(\Sigma)$ is infinitely generated.  Put a
   hyperbolic structure on $\Sigma$ with geodesic boundary.  This
   gives a faithful representation $\sigma \maps \pi_1(\Sigma) \to
   \PSL{2}{\R}$ where the image of the boundary loop $\gamma$ is
   hyperbolic.  As $\pi_1(\Sigma)$ is free, there is a some lift of
   $\til{\sigma}$ to $\til{\PSL{2}{\R}}$.  While
   $\til{\sigma}(\gamma)$ may not have rotation number 0, this time
   $\gamma$ is a primitive element of $H_1(\Sigma)$;  it is easy to
   construct a homomorphism $\til{\rho}$ from $\pi_1(\Sigma)$ to the
   center of $\til{\PSL{2}{\R}}$ so that the rotation number of
   $\til{\rho}(\gamma)$ is minus that of $\til{\sigma}(\gamma)$.  Then
   $\til{\sigma}' \maps g \mapsto \til{\rho}(g) \til{\sigma}(g)$ is
   the required homomorphism.  Note that $\til{\sigma}'$ inherits
   from $\sigma$ the property that only a countable number of points
   have non-trivial stabilizer, and that when the stabilizer is
   non-trivial it is $\Z$.
  
  By viewing $\til{\PSL{2}{\R}}$ as a subgroup of $\homeo(\R) \iso
  \homeo(I)$, we foliate $J$ by $\F$ as a flat foliated bundle with
  holonomy group $\til{\sigma}(\pi_1(\Sigma))$.  As desired, a
  countable set of leaves in the interior of $J$ will be cylinders,
  and all the rest will be planes.  The restriction of the foliation
  $\F$ to the annulus $A$ is a $1$-dimensional foliation $\G$ which is
  a flat foliated bundle over $\gamma = \partial \Sigma$ with holonomy
  generated by $\til{\sigma}(\gamma)$ (see Figure~\ref{fol_of_A}).
  \begin{figure}[h]
    \begin{center}
      \includegraphics[scale=0.5]{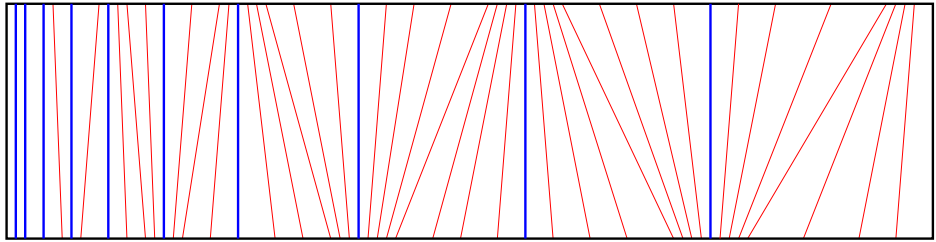}
    \end{center}
    \caption{The element $\til{\sigma}(\gamma)\in \homeo(I)$ has a countable
      collection of fixed points, which suspend to (dark) circles,
      converging to either end.  The holonomy alternates between
      translating in one direction and in the other, which suspends to
      (light) lines which spiral around the (dark) circles. This
      figure depicts the restriction of the foliation $\G$ to one half
      of $A$.}\label{fol_of_A}
  \end{figure}
  The points $t \in I$ parameterize the leaves $\lambda^t$ of $\F$ and
  the leaves $\mu^t$ of $\G$.   Of course, if there is $\beta \in
  \pi_1(\Sigma)$ with $\til{\sigma}(\beta)(t) = s$ then $\lambda^t =
  \lambda^s$, so these labels are not unique. Parameterize $I  = [-1, 1]$
   in such a way that $0$ is an unstable fixed point of
  $\til{\sigma}(\gamma)$.  So the interior of $\lambda^0$ is an
  annulus, and the suspension of the point $0$ in $A$ is a closed loop
  in $\F$ which splits $A$ into two foliated annuli $A^\pm$.
  The interior of each leaf $\lambda^t$ of $\F$ is topologically
  either an annulus or a plane. However, the boundary of each
  $\lambda^t$ consists of the leaves
 \[ 
 \partial \lambda^t = \bigcup_{\beta \in \pi_1(\Sigma)} \mu^{\beta(t)}.
 \]
 For each $t$, the set of boundary components consists of at most a single
 closed loop and countably many lines if $\lambda^t$ is an annulus,
 and countably many lines if $\lambda^t$ is a plane.
 
 Now we'll extend $\F$ slightly into $G$ in the following way.  Let
 $\lambda = \lambda^0$, which is topologically a cylinder, and let
 $\mu^{0\alpha}$ denote the non-circle boundary components of $\lambda$.
 Recall $\mu^0$ denotes the boundary circle of $\lambda^0$.  We split
 $A$ into two foliated annuli $A^\pm$ along the circle $\mu^0$, and
 we denote the foliations by $\G^\pm$.  See Figure~\ref{foliated_A}.
\begin{figure}[h]
  \begin{center}
    \includegraphics[scale=0.4]{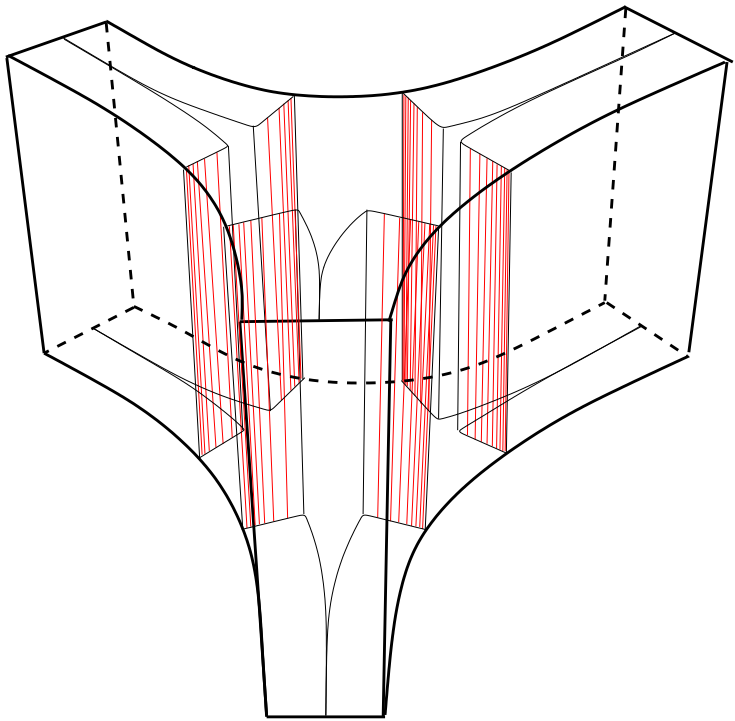}
  \end{center}
\caption{The foliated annuli $A_j$ are split open to annuli foliated by
  $\G_j^\pm$, which are glued together in $G$ by matching their
  holonomy.}\label{foliated_A}
\end{figure}
We extend the foliation to a pair of solid tori $A^+ \times I \subset
G$ and $A^- \times I \subset G$ which are interval bundles over annuli
neighborhoods of $\partial A$ in $\partial G - A$. We foliate these
solid tori $A^\pm \times I$ by the product of the foliations $\G^\pm$
with $I$.  This completes our construction of the foliation of the
interstitial region $J$.
 
 We do this procedure simultaneously for all the
 interstitial regions $J_k$ which bound $G$.  For each object $X$
 constructed in $J$, let $X_k$ denote the corresponding object
 constructed in $J_k$.

 \subsection{Gluing the interstitial foliations} 
 We would like to glue the foliated annuli ends of the $A^\pm_k \times
 I$ together to foliate a product neighborhood of $\partial G -
 \bigcup_k A_k$ (as suggested by Figure~\ref{foliated_A}).  Consider a
 pair of interstitial annuli $A_1$ and $A_2$ which are separated by an
 annulus of $\bdry G - \bigcup_k A_k$.  The task is to glue the
 foliated annulus $A_1^-$ to the foliated annulus $A_2^+$, being
 careful about the topology of the resulting leaves.  Doing the gluing
 amounts to finding a conjugacy between the dynamics of
 $\til{\sigma}_1(\gamma_1)$ on $[-1,0]$ and $\til{\sigma}_2(\gamma_2)$
 on $[0,1]$ with the orientation reversed.  There are many choices of
 such conjugacies; the set of fixed orbits of $\til{\sigma}_1(\gamma_1)$ on
 $[-1,0]$ is a countably infinite discrete sequence accumulating only at 
 $-1$, and the set of fixed orbits of $\til{\sigma}_2(\gamma_2)$ on $[0,1]$
 is a countably infinite discrete sequence accumulating only at $1$. These
 suspend to closed circles in $\G_1^-$ and $\G_2^+$ respectively,
 which can be glued together.  On a complementary interval in
 $A_1^-$ and $A_2^+$, the dynamics of $\til{\sigma}_1(\gamma_1)$ and
 $\til{\sigma}_2(\gamma_2)$ have no fixed points, so the
 set of conjugacies is parameterized by $\homeo(I)$.  We need to
 choose these conjugacies, simultaneously for all adjacent pairs $(A_k^-,
 A_{k+1}^+)$, so that in the resulting foliation the leaves containing
 the middle leaves $\lambda^0_k$ are annuli. 
 
 Consider one of the leaves $\lambda^0_k$, whose boundary lines are
 denoted $\mu^{0\alpha}_k$.  Since all but countably many leaves of
 each $\F_j$ are planes, it is easy to glue so that boundary lines
 $\mu_k^{0\alpha}$ are glued to lines of some $\G_j$ which are
 boundaries of planar leaves $\lambda^{0\alpha}_j$ of $\F_j$.
 However, this is not enough to ensure that the new leaf
 $(\lambda_k^0)'$ which contains $\lambda_k^0$ is an annulus.  The
 leaf $(\lambda_k^0)'$ is built up of pieces which are leaves from the
 $\F_j$, and so far we've only ensured that those pieces that touch
 $\lambda_k^0$ are planar.  Moreover, we have to worry about whether
 the planar pieces of $(\lambda_k^0)'$ could be glued together in a
 cycle, adding to the fundamental group of $(\lambda_k^0)'$.  We'll
 solve these problems by building up the gluings inductively.  More
 precisely, start by gluing one boundary line of $\lambda_k^0$ to a
 planar leaf.  Now glue another one, making sure we don't add any
 fundamental group.  Now glue a boundary component of one of the
 planar leaves we added to a new planar leaf.  Now go work on some
 other $\lambda^0_j$ for a bit, etc.  At each stage, we can avoid
 creating any fundamental group since we've only made finitely many
 gluings.  We can formalize the order of induction as follows.  Note
 that when we're done, each leaf $(\lambda_k^0)'$ will be built up
 from an annulus and infinitely many planes glued together along lines
 in their boundaries.  The dual graph to the pattern of gluings is a
 rooted tree $T_k$, where the root corresponds to the single annulus
 piece, and every vertex has valence the cardinality of the natural
 numbers.  So start with a finite set of such trees $T_k$ and choose a
 bijection from the union of their edges to the natural numbers, so
 that each embedded path from a root vertex outward is labelled by an
 increasing sequence of numbers.  Then do the gluings in the order of
 the labeling.  Since there are only countably many bad choices for
 any given gluing, we can find a good choice arbitrarily close to any
 initial guess, and choose them so that the size of the biggest
 unglued gap goes to zero; thus this process defines the full gluing
 map of all the annuli pairs $(A_k^-, A_{k+1}^+)$.
 
 Next, we'll ``blow air'' into $\F$ along the leaves $(\lambda_k^0)'$ to
 get the final lamination.   We can extend $\F$ to a singular foliation
 $\F^s$ by adding a finite set of annuli $B_j$, parallel to the annuli
 components of $\partial G - \bigcup_k A_k$, which end in pairs on the
 circles $\mu_k^0$. The only singular leaf of $\F^s$ is the single
 branched surface $B$ which is the union of the
 $(\lambda_k^0)'$ and the $B_j$.  Then $B$ is topologically a
 train-track bundle over $S^1$ with fiber the train-track $T$
 obtained from an ideal polygon $P$ by collapsing the cusps of $P$ to
 single branches. It is clear that $B$ can be split open to a
 $\partial P$ bundle over $S^1$, by Denjoying each leaf
 $(\lambda_k^0)'$ before gluing on the annuli $B_j$. Thus $\F^s$ can
 be split open to a nonsingular lamination $\Lambda \subset C$ with a
 single complementary region which is an ideal polygon bundle over
 $S^1$.  Repeating the construction for every complementary region $C$
 of the original $\Lambda$, we obtain a lamination of $M$ all of whose
 complementary regions are ideal polygon bundles over $S^1$.

 \subsection{Tightness is preserved}
 To finish the proof of the Filling Lemma, it remains to check that
 the construction preserves tightness.  We'll make this clear by
 getting an explicit description of the quasi-isometry type of the original
 complementary regions, and understanding how the added leaves are embedded
 in these regions.
 
 Let $\Lambda$ be a genuine lamination with solid torus guts, and take
 a metric on $M$ which makes the leaves of $\Lambda$ hyperbolic.  Let
 $\til{G}$ be the universal cover of a complementary region.  We want to
 describe its quasi-isometry type (with respect to the induced path metric).  Each leaf
 of $\partial \til{G}$ is isometric to $\H^2$. The boundaries of
 interstitial annuli lift to quasigeodesics in leaves of $\partial
 {G}$ bounding quasiconvex regions which are the lifts of the base
 surfaces of the interstitial regions. There are only finitely many
 interstitial annuli, and so the modulus of quasi-isometry of these
 geodesics is uniformly bounded. The gut regions lift to solid tori of
 bounded thickness which are contained in bounded neighborhoods of
 these quasigeodesics in the various boundary leaves of $\partial
 \til{G}$ which they border. Thus, up to coarse quasi-isometry,
 $\til{G}$ is a \emph{tree of hyperbolic planes}. That is, it's
 obtained from countably many copies of $\H^2$,  glued
 together in pairs along convex subsets bounded by complete geodesics.
 A finite number of copies of $\H^2$ come together along a boundary
 geodesic, and the pattern of gluings is treelike---every loop
 in the union is contained in some copy of $\H^2$.
 
 More precisely, each interstitial region is covered by a product
 $\til{\Sigma}_i \times I$ which is uniformly quasi-isometric to a
 region $\til{\Sigma}_i \subset \H^2$ bounded by a collection of
 complete geodesics.  Each leaf $\lambda$ of $\partial \til{G}$ is
 tiled by isometric copies of the $\til{\Sigma}_i$. Pairs of leaves of
 $\partial \til{G}$ are glued together along such tiles, and
 $n$-tuples are glued together around each boundary geodesic of
 $\partial \til{\Sigma}_i$. We can make a graph with two kinds of
 vertices: the first kind correspond to the boundary geodesics of
 tiles $\partial\til{\Sigma}_i$.  The second kind correspond to the
 tiles themselves $\til{\Sigma}_i$.  The edges correspond to a choice
 of a tile and a boundary leaf in that tile.  This gives a graph which
 is a simplicial tree.  (Similar trees of hyperbolic planes arise as
 quasi-isometric models for Cayley graphs of Baumslag-Solitar groups
 \cite{FarbMosher98}, but there the tiles are glued along regions
 bounded by horocycles, not geodesics.)
 
 Looking back at our construction, the final lamination
 $\overline{\Lambda}$ was obtained from $\Lambda$ by inserting new
 leaves into complementary regions to $\Lambda$ so that the added
 leaves $\lambda$ are everywhere transverse to the $I$-bundle regions
 and parallel to the boundary in the gut regions. Such leaves are
 covered in $\til{G}$ by hyperbolic planes which project
 homeomorphically to the $\H^2$ slices in our quasi-isometric model.
 It is clear that a sequence of such slices cannot converge to a pair of
 distinct slices, and so the leaf space of the leaves inside $\til{G}$
 is Hausdorff.  Thus $\overline{\Lambda}$ is tight if $\Lambda$ is.
 This completes the proof of the Filling Lemma.
\halmos \end{proof}

\section{Leaf pocket theorem}

Let $\F$ be a taut foliation with hyperbolic leaves.  Consider a leaf
$\lambda$ of $\F$, and a point $p$ in $\lambda$.  Given a geodesic ray
$r$ starting at $p$, we can ask: Can we choose a transversal $\tau$ at
$p$ so that the holonomy of $\tau$ is defined along all of the ray
$r$?  A yes answer is equivalent to the other leaves not pulling away
from $\lambda$ along $r$ too fast.  Thurston showed the surprising fact that
there are always many directions where the corresponding ray has this
property.  In fact, the set of directions where this is true is dense in
the tangent space to $\lambda$ at $p$.  This section is devoted to
the proof of this theorem, generalized from foliations to
essential laminations.  The existence of so many directions of
``pinching'' of the leaves of $\F$  allows one to piece
together the circles at infinity  of the individual leaves into a
universal circle (see Section~\ref{universal_circles}).

We'll begin with the definition of a \emph{marker}, which is a
way of keeping track of how the leaves of an essential lamination come
together.
\begin{defn}
  Let $\Lambda$ be an essential lamination of $M$ with hyperbolic
  leaves. A {\em marker} for $\Lambda$ is a map
  \[
   m:I \times \R^+ \to \til{M}
  \]
  with the following properties:
\begin{enumerate}
\item{There is a closed set $K \subset I$ such that for each $k \in
    K$, the image of $k \times \R^+$ in $\til{M}$ is a geodesic ray in
    a leaf of $\til{\Lambda}$.  Further, for $k \in I-K$,
    \[
    m(k \times \R^+) \subset \til{M} -\til{\Lambda}.
    \]
    We call these
    rays the {\em horizontal rays} of the marker.}
\item{For each $t \in \R^+$, the interval $m(I \times t)$ is a tight
    transversal.  Further, there is a separation constant $\epsilon$
    for $\Lambda$, such that
    \[
    \mathrm{length}(m(I\times t)) < \epsilon/3.
    \]
    We call these intervals
    the {\em vertical intervals} of the marker.}
\end{enumerate}
\end{defn}

The main theorem of this section is:

\begin{leaf_pocket}\label{leaf_pocket_theorem}
  Let $\Lambda$ be an essential lamination of an atoroidal
  $3$-manifold $M$.  Then for every leaf $\lambda$ of $\til{\Lambda}$,
  the set of endpoints of markers is dense in $S^1_\infty(\lambda)$.
\end{leaf_pocket}

This theorem generalizes Thurston's corresponding theorem for
foliations given in \cite{ThurstonCirclesII}.  The proof here is
completely topological, in contrast to Thurston's proof which relies on
the existence of harmonic measures for foliations.
See Remark~\ref{harmonic} for more on harmonic measures and their relationship to
the Leaf Pocket Theorem.  As such, our conclusion is slightly different from that of
\cite{ThurstonCirclesII}.  Precisely, Thurston showed the
following.  For any $\epsilon>0$ there is a $t>0$ such that if
$\lambda$ and $\mu$ are a pair of leaves which are joined by a
transversal $\tau$ of length $\le t$, then for a random walk $\gamma$
in $\lambda$ starting at $\tau \cap \lambda$, the holonomy transport
of $\tau$ along $\gamma$ has uniformly bounded length with probability
at least $1-\epsilon$.  

The rest of the section is devoted to the proof of the Leaf Pocket
Theorem.
\begin{proof}
  In outline, the proof goes like this.  A \emph{minimal set} of
  $\Lambda$ is a nonempty closed union of leaves which is minimal with
  respect to that property.  A closed union of leaves $\Sigma$ is
  minimal if and only if every leaf in $\Sigma$ is dense in
  $\Sigma$.  The main case of the theorem is for leaves of
  $\til{\Lambda}$ which cover leaves in a minimal set.  While not
  every leaf of $\Lambda$ is contained in a minimal set, the closure
  of every leaf in $\Lambda$ contains some minimal set.  Once we know
  the theorem for leaves in minimal sets, it is not hard to prove the
  theorem for all leaves.
  
  To continue the outline, let $\Sigma$ be a minimal set of $\Lambda$.
  First we  show that some leaf of $\Sigma$ has nontrivial
  fundamental group.  Then there is a closed geodesic $\gamma$ in a
  leaf of $\Sigma$.  Using the holonomy around $\gamma$, we construct
  an ``immersed sawblade'' in $\Lambda$, from which we build a
  collection of markers in $\til{\Lambda}$.  We show that any leaf in
  $\til{\Sigma}$ intersects at least one of these markers.  We then 
  show that in any given leaf $\lambda$ of $\til{\Lambda}$, the
  endpoints of these markers is dense in $S^1_\infty(\lambda)$.
    
  Let's begin the proof in the case of a minimal set $\Sigma$ of
  $\Lambda$, that is, for leaves $\lambda$ of $\til{\Lambda}$ which lie
  in the inverse image $\til{\Sigma} \subset \til{\Lambda}$ of
  $\Sigma$.

  \subsection{Each leaf of $\til{\Sigma}$ has a marker}
  
  First, we claim that $\Sigma$ contains a leaf with a nontrivial
  closed geodesic.  Suppose to the contrary that every leaf of
  $\Sigma$ is simply connected.  Since $M$ is atoroidal, the
  lamination $\Lambda$ has hyperbolic leaves, and so each leaf of
  $\Sigma$ is a hyperbolic plane.  Since the holonomy on $\Sigma$ is
  trivial, it's easy to construct a nontrivial invariant transverse
  measure supported on $\Sigma$.  In codimension $1$, Plante showed
  that leaves in the support of an invariant transverse measure have
  polynomial area growth \cite{Plante75}.  This is a contradiction as
  the hyperbolic plane has exponential growth.  So there is a leaf
  $\lambda$ of $\Sigma$ which contains a closed geodesic $\gamma$.

  (Note: In \cite{Plante75}, Plante was concerned only with foliations,
  but his proof applies to laminations as well.  The only new issue is
  the step which involves passing to a \emph{co-orientable} finite cover,
  since not all laminations are virtually co-orientable.  However, one
  can find an open neighborhood $N(\Lambda)$ of $\Lambda$ where
  $\Lambda$ can be extended to a foliation, and pass if necessary to a
  finite cover of $N(\Lambda)$ where the foliation becomes
  co-orientable. That is, one can find a $3$-manifold $\widehat{N}$
  and a compact co-oriented lamination $\widehat{\Lambda}$ in
  $\widehat{N}$ which maps by either an embedding or a double cover to
  a submanifold of $M$, taking $\widehat{\Lambda}$ to $\Lambda$.
  Plante's theorem applies to $\widehat{\Lambda}$, and therefore also
  to $\Lambda$.)
  \begin{figure}[h]
    \begin{center}
      \includegraphics[scale=0.4]{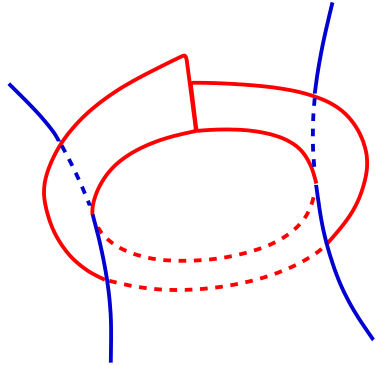}
    \end{center}
    \caption{Holonomy transport of some small transversal around $\gamma$ sweeps out an
      immersed sawblade.}\label{sawblade}
  \end{figure}
  
  Now we'll use $\gamma$ to show that each leaf of $\til{\Sigma}$
  intersects a marker.  We think of $\gamma$ as an interval, whose
  endpoints are the same point in $\lambda$.  If $\Sigma$ is a closed
  isolated leaf then constructing markers is trivial, so we assume this
  is not the case.  As $\lambda$ is in a minimal set, the geodesic
  $\gamma$ is a limit of leaves of $\Lambda$ on at least one side.
  Let $\tau$ be a tight transversal with lower endpoint in $\lambda$
  to a non-isolated side of $\lambda$, sufficiently small so that the
  holonomy transport of $\tau$ around $\gamma$ exists.  Moreover, if
  $\lambda$ is not isolated in $\Sigma$, choose $\tau$ on a non-isolated
  side in $\Sigma$.  Then, by minimality of $\Sigma$, every leaf of
  $\Sigma$ intersects $\tau$.
  
  Consider the holonomy transport $T\maps I \times \gamma \to M$ of $\tau$
  around $\gamma$.  Even though $\gamma(0)=\gamma(1)$, we do not
  usually have $T(I,\gamma(0)) = T(I,\gamma(1))$. However, either
  $T(I,\gamma(1)) \subset T(I,\gamma(0))$ or vice versa.  See
  Figure~\ref{sawblade}.  By reversing the orientation on $\gamma$ if
  necessary, we can ensure $T(I,\gamma(1)) \subset T(I,\gamma(0))$.
  We will call the image of $T$ an {\em immersed sawblade}, and denote
  it $G$.  The universal cover $\til{G}$ of $G$ lifts to $\til{M}$
  where it runs along a lift of $\lambda$.
  \begin{figure}[h]
    \begin{center}
      \includegraphics[scale=0.4]{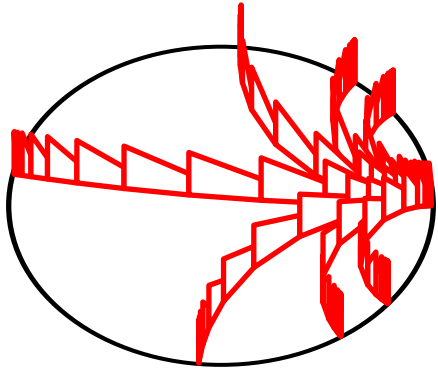}
    \end{center}
    \caption{Lifts of a sawblade $\til{G}$ running along $\lambda$.}
    \label{sawblade2}
  \end{figure}
  See Figure~\ref{sawblade2} for one possible configuration which
  makes clear why $G$ is called an immersed sawblade.  If we take
  $\tau$ short enough so that $G$ has height less than $1/3$ a
  separation constant for $\Lambda$, a lift of the immersed sawblade
  $G$ must contains the image of a marker $m$.  More precisely,
  consider any lift of $\tau$ to $\til{M}$ and look at the lift of
  $\til{G}$ which contains it.  The sawblade structure of $\til{G}$
  forces holonomy of $\tau$ to be defined for all time in the positive
  direction along the lift of $\gamma$.  Sweeping the lift of $\tau$
  along the lift of $\gamma$ gives a marker $m$ which is contained in
  $\til{G}$.  Actually, $m$ isn't quite a marker because the
  definition requires that each horizontal ray be a geodesic.
  Shrinking $\tau$ and thus $G$ if necessary, we can ensure that each
  horizontal ray has uniformly small geodesic curvature in its leaf.
  Then each horizontal ray of $m$ is a $K$-quasigeodesic for a uniform
  $K$ which we can make arbitrarily close to $1$.  We can then
  straighten these to genuine geodesics to construct a real marker.
  Making $K$ small enough, the final marker lies in a small
  $\delta$-neighborhood of $\til{G}$.

  Taking the possible lifts of $\tau$ to $\til{M}$ gives us an
  equivariant family $\mathscr{M}$ of markers in $\til{M}$.  Now
  consider any $\lambda'$ in $\Sigma$ and some leaf $\til{\lambda}'$
  in $\til{\Sigma}$ which covers it.  As noted above, $\lambda'$ must
  intersect the transversal $\tau$.  Thus one of the markers in
  $\mathscr{M}$ intersects $\til{\lambda}'$.  So every leaf of
  $\til{\Sigma}$ intersects a marker.

  \subsection{Endpoints of markers}  
  
  Consider the family $\mathscr{M}$ of markers coming from the lifts
  of $\tau$ above.  Let $\lambda$ be any leaf of $\til{\Sigma}$.
  Consider the collection of geodesic rays in $\lambda$ which are
  horizontal rays of markers of $\mathscr{M}$.  These rays originate
  from the intersections of lifts of $\tau$ with $\lambda$.  We want
  to show that the endpoints of these rays are dense in
  $S^1_\infty(\lambda)$.  Equivalently, let $\mathscr{G}$ denote the
  set of all lifts of $\til{G}$ to $\til{M}$, and look at the
  intersections of $\mathscr{G}$ with $\lambda$.  We need to show that
  the endpoints of these quasigeodesic rays are dense in
  $S^1_\infty(\lambda)$.
  
  First, we'll show there are infinitely many intersections of
  $\mathscr{G}$ with $\lambda$.  Since $\Sigma$ is minimal, there is a
  uniform $R$ such that for any point $p$ in any leaf of $\Sigma$, the disk
  in that leaf about $p$ of radius $R$ intersects $\tau$.  This is
  because the function of leafwise distance to $\tau$ is upper
  semi-continuous on $\Sigma$.  Applying this to $\lambda$, we have
  that any point $p \in \lambda$ is within $R$ of a lift of $\tau$,
  and therefore within $R$ of an intersection of $\mathscr{G}$ and
  $\lambda$.  Thus, there are infinitely many quasigeodesic rays
  $\gamma_i$ contained in $\lambda \cap \mathscr{G}$.
  
  Next, we'll show that any two distinct $\gamma_i$ and $\gamma_j$
  have disjoint endpoints in $S^1_\infty(\lambda)$ (in particular,
  Figure~\ref{sawblade2} depicts an impossible configuration). 
  Downstairs in $M$, choose a $\delta$ so
  that the leafwise $\delta$-neighborhood of $G$ deformation retracts
  to $G$.  Note that $\gamma_i$ and $\gamma_j$ come from distinct
  lifts of $\til{G}$ in $\mathscr{G}$, as otherwise some lift of
  $\tau$ intersects $\lambda$ twice.  Lifting the
  $\delta$-neighborhood of $G$ upstairs, we see that in $\lambda$ the
  $\delta$-neighborhoods of $\gamma_i$ and $\gamma_j$ are disjoint.
  By shrinking $\tau$ if necessary, we can make the uniform
  quasigeodesic constant for all the $\gamma_k$ so small, that the
  $\delta/3$ neighborhood of any $\gamma_k$ contains a geodesic ray.
  If $\gamma_i$ and $\gamma_j$ had a common endpoint, consider the
  corresponding geodesic rays that they are $\delta/3$-close to.
  There are points on these geodesic rays that are less than
  $\delta/3$ apart, and so the distance between $\gamma_i$ and
  $\gamma_j$ is less than $\delta$.  But this contradicts that
  $\gamma_i$ and $\gamma_j$ have disjoint $\delta$-neighborhoods in
  $\lambda$.  So $\gamma_i$ and $\gamma_j$ have disjoint endpoints.
  
  From now on, we'll replace the original $\gamma_i$ by the
  corresponding geodesics.  Note that the $\gamma_i$ have disjoint
  $\delta/3$-neighborhoods, and that the $\gamma_i$ are an $R +
  \delta$ net for $\lambda$.  Now suppose that the endpoints of the
  $\gamma_i$ are not dense.  Let $I$ be an open interval in
  $S^1_\infty(\lambda)$ containing no endpoints.  Let $J$ be a closed
  interval contained in $I$.  Let $H$ be the half-space associated to
  the smaller interval $J$.  As the $\gamma_i$ are a net, there must
  be an infinite sequence of disjoint $\gamma_i$, call them
  $\alpha_k$, which intersect $H$.  Passing to a subsequence, we can
  assume that the terminal endpoints of the $\alpha_k$ converge to a
  point $p$ in $S^1_\infty(\lambda)$, and that the initial endpoints
  converge to a point $q$ in $\lambda \cup S^1_\infty(\lambda)$.  The
  point $q$ is in $H \cup J$, and, by assumption, $p$ can't be in $I$.
  So $p \neq q$.  But then the $\alpha_k$ converge to the geodesic ray
  joining $q$ to $p$, which contradicts that the $\alpha_k$ have
  disjoint $\delta/3$-neighborhoods.  Thus the endpoints of markers
  are dense in $S^1_\infty(\lambda)$.  This proves the theorem in the
  case of leaves in minimal sets.
  
  \subsection{Properties of markers}
  Before moving on, we note that we actually know more than just that
  the endpoints of the markers are dense in $S^1_\infty(\lambda)$.  If $\Sigma$
  is not a compact leaf, then every leaf of $\Sigma$ intersects the
  \emph{interior} of $\tau$ many times.  Thus there must be a marker
  $m$ for each leaf $\lambda$ in $\til{\Sigma}$ so that $\lambda$ is
  in the interior of marker (in the sense that the vertical intervals
  of the marker intersect $\lambda$ in their interior).  If, on the
  other hand, $\Sigma$ is a compact leaf then all the markers we
  construct have endpoints in $\Sigma$, and extend only to one side of
  $\Sigma$.  In this case, we repeat the construction on the other
  side of $\gamma$ to get dense sets of markers extending in both
  directions.

  \subsection{Leaves not in minimal sets}
  The basic idea here is that leaves that are not in minimal sets get
  very close to leaves which are in minimal sets, and pick up markers
  that way.  Fix a leaf $\lambda$ of $\Lambda$, and a point $p$ in
  $\lambda$.  For each unit vector $v$ tangent to $\lambda$ at $p$, we
  need to show that there is a marker intersecting $\lambda$ whose
  endpoint in $S^1_\infty(\lambda)$ is arbitrarily close (in the
  visual sense) to $v$.  Fix an $\epsilon > 0$.  Let $S_\epsilon$
  be the sector of $\lambda$ bounded by two geodesic rays starting at
  $p$, so that the angle between the two rays is $\epsilon$, and so
  that $v$ points down the center of $S_\epsilon$.  We claim the closure of
  $S_\epsilon$ in $\Lambda$ must contain a minimal set.  Pick
  a sequence of points $m_i$ in $S_\epsilon$ marching out towards
  infinity so that the disk about $m_i$ in $\lambda$ of radius $i$ is
  contained in $S_\epsilon$.  If $m$ is any limit point in $M$ of the
  $m_i$, the closure of $S_\epsilon$ contains the entire leaf
  containing $m$, and thus contains any minimal set in the
  closure of that leaf.   Let $\Sigma$ be such a minimal set.

  Fix a point $q$ in a leaf $\sigma$ of
  $\Sigma$.  As $q$ is in the closure of $S_\epsilon$, there exist points $p_i$
  in $S_\epsilon$ which converge to $q$ in $M$.  We can assume that locally
  near $q$, the $p_i$ all lie to a fixed side of $\sigma$.  Let $v_i$
  be the unit tangent vector to $\lambda$ at $p_i$ which points
  directly away from the base point $p$.  By passing to a subsequence,
  we can assume that the $v_i$ converge to a vector $w$ tangent to
  $\sigma$ at $q$.  As $\sigma$ is in a minimal set, there
  exists a marker whose endpoint in $S^1_\infty(\sigma)$ makes an
  angle of $< \epsilon/2$ with $w$.  Initially, that marker might not
  come near the point $q$, but at the cost of shortening the vertical
  intervals we can drag it over so that one of the vertical intervals
  contains $q$.  Moreover, we can choose this marker so that it lies
  to the same side of $\sigma$ as the $p_i$'s.  This marker intersects
  $\lambda$ in horizontal rays which nearly contain the geodesic rays
  given by $(p_i, v_i)$.  For $i$ large enough, these horizontal rays
  must be contained in the enlarged sector $S_{2\epsilon}$.  Thus the
  endpoints of these markers make a visual angle of at most $2
  \epsilon$ with our initial vector $v$.  It follows that the
  endpoints of markers are dense in $S^1_\infty(\lambda)$.
\halmos \end{proof}

\begin{rmk}\label{super_density}
  Note that the proof actually shows a little more about the set of
  markers than just that their endpoints are dense in
  $S^1_\infty(\lambda)$.  Namely, if $\lambda$ is non-compact, then
  the endpoints of markers which intersect $\lambda$ in
  their interior are dense.  If $\lambda$ is compact, then the
  endpoints of markers which extend to a fixed side of $\lambda$ are
  dense.
\end{rmk}

\begin{rmk}
  Notice that the ``essential'' hypothesis in this theorem is
  excessive. All that is necessary is that $\Lambda$ is compact,
  codimension one, and has hyperbolic leaves.  Also, while
  Theorem~\ref{leaf_pocket_theorem} generalizes Thurston's original
  theorem from foliations to all essential laminations, this does not
  allow one to generalize his theorem on the existence of universal
  circles to the case of genuine laminations.
\end{rmk}

\begin{rmk}\label{harmonic}
  A more analytic proof of this theorem might run along the following
  lines: The set of harmonic probability measures supported on a
  compact lamination $\Lambda$ is a compact convex set; the extremal
  points in this set (those which can't be written as a nontrivial
  convex sum of harmonic probability measures) are the ergodic
  harmonic measures \cite{Garnett83,Candel:Garnett}.  The holonomy
  transport along a random walk in a leaf preserves the
  (infinitesimal) harmonic measure of a transversal, on average.  It
  follows by an analytic argument that the holonomy of some
  transversal does not blow up for a.e.~random walk on a leaf in the
  support of an ergodic harmonic measure.  If $\lambda$ is a leaf not
  in the support of any ergodic harmonic measure, a random walk on
  $\lambda$ is dispersive, and must eventually wander arbitrarily
  close to the support of an ergodic measure.
\end{rmk}

\section{Universal circles}\label{universal_circles}

Let $M$ be an orientable 3-manifold with a taut foliation $\F$.  In
the last section, we saw that the leaves of $\til{\F}$ come together
in many directions.  In this section, we'll  explain how this allows the
construction of a master \emph{universal circle} by assembling the
circles at infinity of the leaves of $\til{\F}$.  This universal
circle comes equipped with an action of $\pi_1(M)$.  A simple example
to keep in mind is that of a surface bundle over $S^1$, where one can
identify the circles at infinity of any two leaves using the product
structure of $\til{M}$.  In this case, the universal circle can be
identified with any particular circle at infinity.  In general, the
branching behavior of the leaf space makes the relation between the
universal circle and the circle at infinity of a leaf more
complicated.

We'll start with the precise definition of a \emph{universal circle}.
Recall that two leaves are in $\til{\F}$ are \emph{comparable} if they can be
joined by an arc transverse to $\til{\F}$.  Then:
\begin{defn} \label{univ_circle_def}
A {\em universal circle}, $S^1_\un$, for $\F$ is the following data:
\begin{enumerate}
\item A representation 
  \[
  \rho_\un \maps \pi_1(M) \to \homeo(S^1_\un).
  \]
\item For every leaf $\lambda$ of $\til{\F}$, a monotone map
\[
  \phi_\lambda \maps S^1_\un \to S^1_\infty(\lambda).
\]

\item For every $\alpha \in \pi_1(M)$, the following diagram commutes:
\[
\xymatrix{S^1_\un   \ar[d]_{\phi_{\lambda}}  \ar[r]^{\rho_\un(\alpha)}  & 
  S^1_\un \ar[d]_{\phi_{\alpha(\lambda)}}\\
  S^1_\infty(\lambda)    \ar[r]^{\alpha}   &     S^1_\infty(\alpha(\lambda)) }
\]
\label{diagram_commutes}

\item For any leaf $\mu$, the \emph{associated gaps} are the maximal
  connected intervals in $S^1_\un$ mapped to points by $\phi_\mu$.
  The complement of the gaps in $S^1_\un$ is the {\em core} associated
  to $\mu$.  We require that for any pair $(\mu,\lambda)$ of
  incomparable leaves, the core associated to $\lambda$ is contained
  in a single gap associated to $\mu$, and vice versa.  \label{assoc_gaps}
\end{enumerate}
\end{defn}

The purpose of this section is to prove:

\begin{Theorem}[Thurston]\label{universal_circle_exists}
  Let $\F$ be a taut foliation of an orientable $3$-manifold with
  hyperbolic leaves. Then there exists a universal circle for $\F$.
\end{Theorem}

Moreover, one has:

\begin{Theorem}\label{universal_circle_action_faithful}
  If $M$ is atoroidal, the action of $\pi_1(M)$ on the universal
  circle is faithful.
\end{Theorem}

Theorem~\ref{universal_circle_exists} is announced in
\cite{ThurstonCirclesII}, but unfortunately this paper is mostly
unwritten. Thurston has outlined a construction of $S^1_\un$ in
several lectures, and given many details. An alternate construction of
$S^1_\un$ when $\F$ is co-orientable is given in
\cite{Calegari:thesis}.  As a public service, we give a complete proof
here.  The basic idea is to look at certain sections of the circle bundle
$E_\infty$ over the leaf space $L$, namely the ``special sections''
which are canonically defined.  The set of all special sections,
$\mathscr{S}$, has a natural circular order, and can be completed to
form $S^1_\un$.  The map $S^1_\un \to S^1_\infty(\lambda)$ comes from
restricting sections in $\mathscr{S}$ to $\lambda$.  Since special
sections are canonically defined, they are invariant under the action
of $\pi_1(M)$ on $E_\infty$, and this gives the action of $\pi_1(M)$
on $S^1_\un$.  The hard work in proving
Theorem~\ref{universal_circle_exists} is all in the Leaf Pocket
Theorem~\ref{leaf_pocket_theorem}, which is used to define the
correct sections.  The arguments in this section are mostly formal.

\subsection{Monotone maps of circles}
We'll warm up to our study of special sections by investigating some generalities
about monotone maps between circles. 
A {\em monotone map} $\phi \maps S^1 \to S^1$ is a degree one map
which does not reverse (but might degenerate) the cyclic ordering on
triples of points.  Note that it is implicit in this definition that we are
considering maps between {\em oriented} circles.
A map is monotone if for every $p \in S^1$, the
preimage $\phi^{-1}(p)$ is contractible. This use of the word monotone agrees
with the standard use from decomposition theory, in the sense of R. L. Moore.

\begin{defn}
A {\em monotone relation} between an ordered pair of circles $S^1_1$ and $S^1_2$ is a third
circle $S^1_{12}$ and two monotone maps $\phi_i \maps S^1_{12} \to S^1_i$
for $i = 1,2$.
\end{defn}
For ease of notation, we denote a monotone relation between two
circles by the name of the source of the two monotone maps to these
circles. The monotone maps are part of the data, of course.

The following lemma gives a pushout construction for monotone maps.

\begin{Lemma}\label{composition}
  Consider three circles $S^1_1, S^1_2$, and $S^1_3$.  Let $S^1_{12}$
  be a monotone relation between $S^1_1$ and $S^1_2$, and $S^1_{23}$ a
  relation between $S^1_2$ and $S^1_3$.  Then there is a canonical
  (leftmost) monotone relation $S^1_{13}$ between $S^1_1$ and $S^1_3$.
\end{Lemma}
\begin{proof}
  The mapping cylinder of a monotone relation is literally a cylinder
  whose interior is foliated as a product, but for which distinct
  intervals of leaves converge to a single point in the target circle
  (see Figure~\ref{map_cylinder}).
\begin{figure}[h]
  \begin{center}
    \includegraphics[scale=0.8]{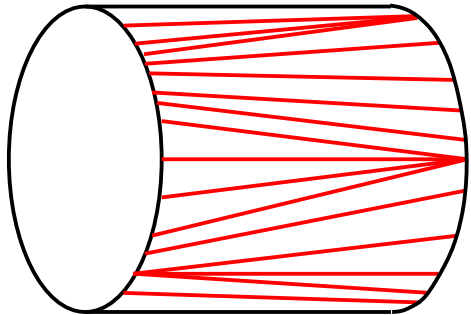}
  \end{center}
  \caption{The mapping cylinder of a monotone relation.}\label{map_cylinder}
\end{figure}
The mapping cylinders of our two monotone relations can be glued along
the intermediate $S^1_2$ to give a cylinder $C$ with a foliation which
is singular along the middle and two boundary circles.  Using a
leftmost rule, we will resolve these singularities to construct the
mapping cylinder of the composite.
\begin{figure}[h]
  \begin{center}
    \includegraphics[scale=0.5]{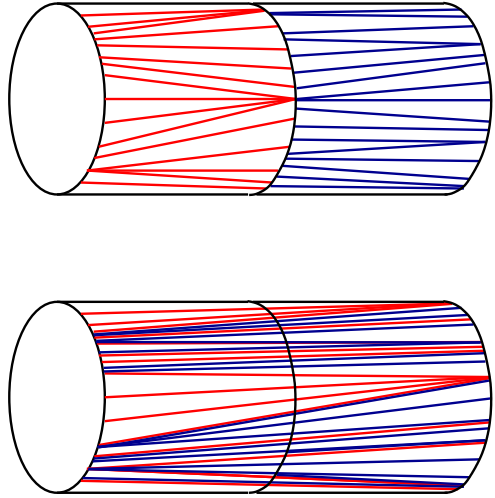}
  \end{center}
\caption{Monotone relations can be composed by splitting the
associated singular foliations according to the {\em leftmost rule}.}\label{compose_monotone}
\end{figure}
If $p \in S^1_2$ is singular, the leaf through $p$ consists of the
union of a cone and a line, or of two cones.  (Throughout, consult
Figure~\ref{compose_monotone}.)  The boundary of this singular leaf is
a train track with $3$ or $4$ branches and a single switch at $p$.  If
the train track has $3$ branches, we can split it open until the
switch is on a boundary circle of $C$. If the track has $4$ branches,
we resolve the switch by pushing the upper left branch over the right
branch, and then split open as before. The regions bounded by split
open tracks can be foliated by a product foliation. The result is a
foliation of the interior of $C$ which is singular along the two
boundary circles, which are canonically identified with $S^1_1$ and
$S^1_3$.  This gives the needed relation.
\halmos \end{proof}

If we think of monotone relations between \emph{ordered pairs}, then
Lemma~\ref{composition} shows us how to compose monotone relations.
The topological realization of the pushout as a canonical splitting of
a singular foliation of a cylinder implies that the composition of
monotone relations is {\em associative}. 

Lemma~\ref{composition} generalizes to continuous families of monotone
maps, as follows.  The reader may wish to skip ahead to
Section~\ref{leftmost_section} to see how this material is needed in the sequel
before proceeding.  Consider a continuous family of monotone maps 
\[
  \phi_t \maps S^1\times t \to S^1 \times t
\]
parameterized by $t$ in an interval $I$.  Let $C=S^1 \times I$ be a
parameterized cylinder, and let $\phi \maps C \to C$ be the map associated
to the $\phi_t$, which preserves the $I$-coordinate.

The cylinder $C$ is foliated by $\{\text{point}\} \times I$.  The
image of this foliation under $\phi$ is denoted by $\F_\phi$, and
$\F_\phi$ is said to be \emph{monotonely parameterized by $\phi$}. The
foliation is singular, since distinct leaves are not disjoint.  While
leaves may overlap, they can't cross.  It is clear from the definition
that a monotone parameterization of $C$ gives a monotone relation
between the circles of $\partial C$.   
\begin{Lemma}\label{monotone_parameterization}\label{monotone_param}
  Let $C = S^1 \times I$ be a parameterized cylinder, and let $\F$ be
  a singular foliation of $C$ transverse to every $S^1 \times
  \{\text{point}\}$. Then $\F$ is monotonely parameterized by a
  canonical leftmost $\phi$.
\end{Lemma}
\begin{proof}
  We can think of the singular foliation $\mathscr{F}$ as a foliation by
  train tracks, in the generalized sense that continuous families of
  branches may coalesce at a switch. A dense subset of $\F$ can be
  exhausted by finite train-tracks $F_i$.  These can be
  desingularized, pushing left branches over right branches to the
  right boundary where they can be split open, and right branches
  under left branches to the left boundary where they can be split
  open. This gives a family of trivial (product) laminations in $C$ of
  the form $S_i \times I$ where $S_i$ is a finite set, which maps by a
  monotone family of maps to the train track $F_i$. As the $F_i$
  accumulate in $C$, so do the $S_i$, and in the limit 
  we construct a monotone family of maps from $C$ foliated as a
  product to $C$ foliated by $\F$.  This limiting map gives our
  canonical monotone parameterization of $\F$.  It is easy to check
  that the map does not depend up to homeomorphism on the choice of
  exhaustion by $F_i$, and is natural.
\halmos \end{proof}

\subsection{Leftmost sections for comparable leaves}\label{leftmost_section}

Fix a co-orientation of $\til{\F}$.  For a pair $(\lambda, \mu)$ of
comparable leaves in $\til{\F}$, we will write $\lambda > \mu$ to mean
that $\lambda$ lies above $\mu$ with respect to our co-orientation.
Consider such a pair $\lambda > \mu$.  Let $I \subset L$ be the
interval joining $\lambda$ to $\mu$. Then $C=E_\infty|_I$ is
topologically a cylinder foliated by circles $S^1_\infty(\nu)$ for
$\lambda \geq \nu \geq \mu$.  In this subsection, we will define
certain sections of $C$, the \emph{leftmost sections}. The leftmost
sections of $C$ are defined in terms of the endpoints of markers, so
we'll begin with a discussion of markers.  Let $m$ be a marker for
$\til{\F}$ which intersects only leaves in $I$.  Then the endpoints of
the horizontal rays of $m$ form an interval in $C \subset E_\infty$
which is transverse to the circle fibers.  From now on, we will abuse
notation and refer to the interval of endpoints in $E_\infty$ as a
marker.  By the Leaf Pocket Theorem~\ref{leaf_pocket_theorem}, the
markers are dense in each circle fiber $S^1_\infty(\nu)$ of $C$.
Moreover, even if we only look at markers which extend to a fixed side
of $S^1_\infty(\nu)$, then these markers are dense in
$S^1_\infty(\nu)$ (see Remark~\ref{super_density}).  The key to
understanding the whole set of markers is the following:
\begin{Lemma}\label{marker_crossing}
  Let $m_1$ and $m_2$ be two markers in $C$.  If $m_1$ and $m_2$ are
  not disjoint, their union $m_1 \cup m_2$ is also an interval transverse
  to the circle fibers (in particular, $m_1 \cap m_2$ is a subinterval of each $m_i$).  
\end{Lemma}

\begin{proof}
  Suppose the markers intersect.  Let $\lambda$ be a leaf where the
  endpoints of $m_1$ and $m_2$ are the same, and $\mu$ be a leaf where
  the endpoints of the $m_i$ differ.  Remember that we required that
  the vertical intervals in our markers are very thin, and have length
  $< \epsilon/3$ for some separation constant $\epsilon$.  In
  particular, the vertical intervals of the $m_i$ from $\lambda$ to
  $\mu$ have length $< \epsilon/3$.  Let $\gamma^\lambda_i$ denote the
  horizontal ray of $m_i$ in $\lambda$ and similarly with
  $\gamma^\mu_i$.  Shortening the markers horizontally if necessary,
  we can make $\gamma^\lambda_1$ lie in an $\epsilon/3$ neighborhood
  of $\gamma^\lambda_2$.  If we start at some point $p$ on
  $\gamma_1^\mu$, we can join it to a point $q$ in $\gamma_2^\mu$ by a
  path of length $< \epsilon$ by going up along $m_1$ to
  $\gamma_1^\lambda$, over to $\gamma_2^\lambda$ in $\lambda$, and
  down $m_2$ to $\gamma_2^\mu$.  But the $\gamma^\mu_i$ diverge in
  $\mu$, so the distance in $\mu$ between $p$ and $q$ can be made
  arbitrarily large.  But the distance in $\til{M}$ between $p$ and
  $q$ is at most $\epsilon$.  This contradicts that $\epsilon$ is a
  separation constant for $\F$, as $\mu$ is not
  quasi-isometrically embedded in its $\epsilon$-neighborhood.
\halmos \end{proof}

By the above lemma, we can amalgamate intersecting markers into one
larger interval.  If we take a maximal such union of markers, we get
an interval $m$ which is quite possibly open at either end.  We'll
abuse notation and call such an $m$ a marker.  Because of the density of
markers, the ends $m$ can't wiggle violently, and the closure of $m$
is a closed interval transverse to the circle fibers of $C$.  
The closures of two maximal markers can intersect only in at most their
endpoints.  Thus, the set of all markers gives us something which
approximates a singular foliation of $C$, not unlike
Figure~\ref{compose_monotone}.  Later, we'll describe how to
``integrate'' the set of markers into a proper singular foliation.

Now we're in a position to define the \emph{leftmost sections} of $C$.
First, an \emph{admissible section} is a section $\tau \maps  I \to C$ whose
image does not cross, but might run into, any marker.  The
\emph{leftmost section} starting at $p \in S^1_\infty(\mu)$ is an
admissible section $\tau \maps I \to C$ which is anticlockwisemost
among all such sections in the following sense: for any leaf $\nu \in
I$, the value $\tau(\nu)$ is anticlockwisemost among all endpoints of
admissible sections starting at $p$ and ending on $S^1_\infty(\nu)$.
\begin{Lemma} 
  The leftmost section exists and is continuous.
\end{Lemma}
\begin{proof}
  For a fixed $p$ in $S^1_\infty(\mu)$, we will construct a series of
  approximations to the leftmost section.  Consider a finite set of
  markers $\mathscr{M}$ which intersects every leaf $S^1_\infty(\nu)$.
  Consider paths $\gamma \maps [0,1] \to C$ which are almost sections
  in the sense that the composition of $\gamma$ and the projection
  $[0,1] \to I$ is a monotone surjection.  It makes sense to say such a
  path is admissible.  Define $\tau_\mathscr{M}$ to be the leftmost
  admissible path in this sense.  Explicitly, $\tau_\mathscr{M}$ is
  constructed by starting at $p$, heading left until we hit a marker
  $m \in \mathscr{M}$, going up along $m$ until it ends, then left
  again until we hit a marker, etc.   See Figure~\ref{appox_leftward}
  \begin{figure}[htb]
    \psfrag{p}{$p$}
    \psfrag{C}{$C$}
    \begin{center}
      \includegraphics[scale=0.5]{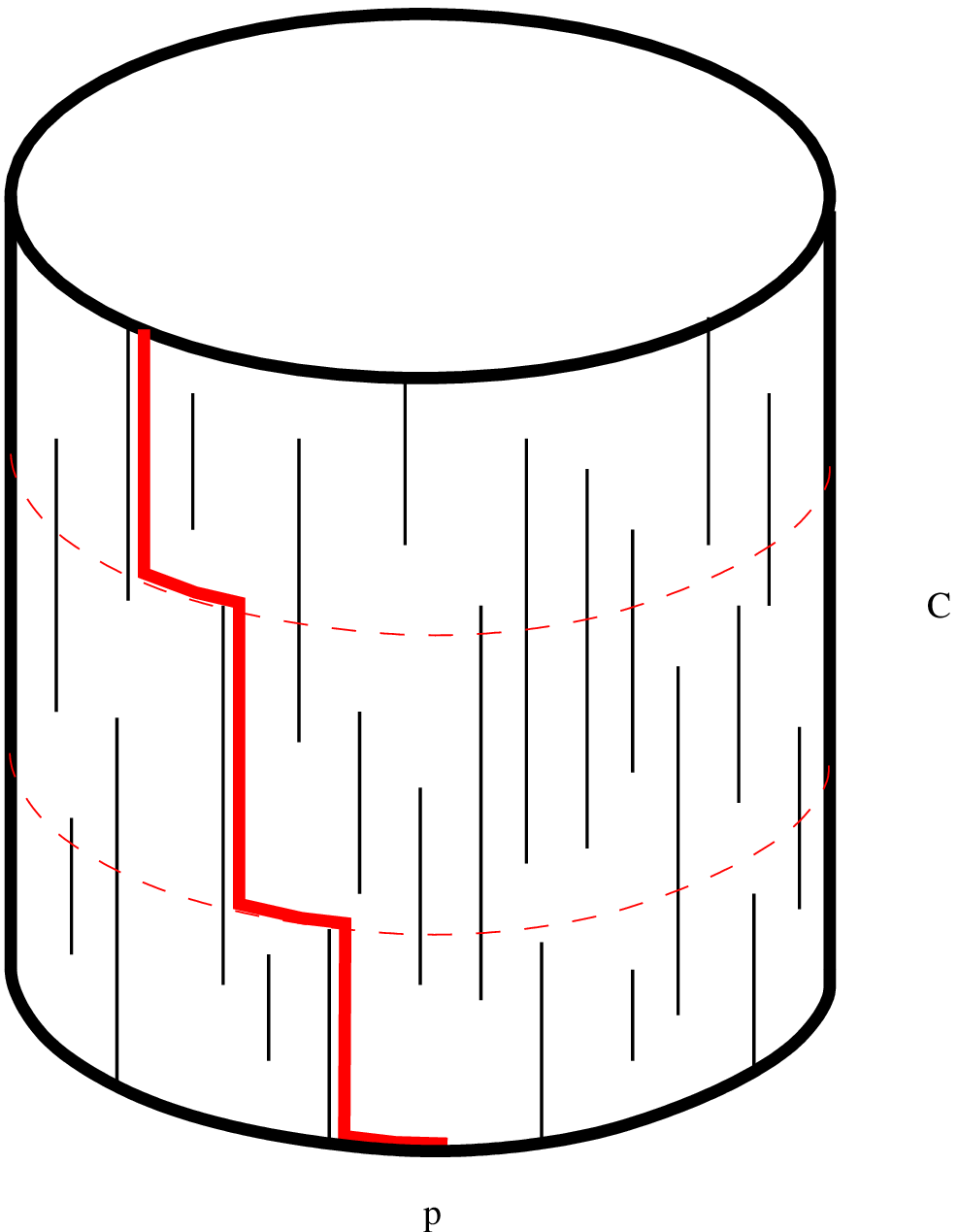}
    \end{center}
    \caption{An approximate leftmost section $\tau_\mathscr{M}$.}
    \label{appox_leftward}
  \end{figure}
  
  Note that if we add additional elements to $\mathscr{M}$, the
  intersection of $\tau_\mathscr{M}$ with any fixed $S^1_\infty(\nu)$
  moves to the right.  The leftmost section $\tau$ to the full set of
  markers is essentially the righthanded limit of all the
  $\tau_\mathscr{M}$.  To be precise, let's work in the universal
  cover of $C$ which is $\R \times I$, and consider the
  $\tau_\mathscr{M}$ to be based at some fixed lift $\til{p}$ of $p$.
  Then we define a section $\tau \maps I \to C$ by
  \[
  \tau(\nu) = \sup_\mathscr{M} \big(  \min( \tau_\mathscr{M} \cap S^1_\infty(\nu)) \big).
  \]
  Here's why this supremum exists: First note that we can restrict the
  supremum to only markers $\mathscr{M}$ which contain some fixed
  set of markers $\mathscr{M}_0$.  Then, for any $\mathscr{M} \supset
  \mathscr{M}_0$, the path $\tau_\mathscr{M}$ lies to the left of the
  the rightward analogue of $\tau_{\mathscr{M}_0}$.
  
  Because of the density of markers in {\em every} circle $S^1_\infty(\lambda)$,
  it is not hard to see that
  $\tau$ is continuous.  Since the supremum was taken over all finite
  subsets $\mathscr{M}$, the section $\tau$ is admissible with respect
  to the full set of markers.  Finally, if $\gamma$ is any admissible
  section, it lies to the right of any $\tau_\mathscr{M}$, and hence
  to the right of $\tau$.  So $\tau$ is the promised leftmost section
  starting at $p$.
\halmos \end{proof}

While a given leftmost section is continuous, the leftmost sections
may not vary continuously as a function of $p$.   However,  any two
leftmost sections do not cross, though they can coalesce.
Symmetrically, given a point $q$ in the top boundary component
$S^1_\infty(\lambda)$, we can talk about the rightmost (downward)
section of $C$ starting at $q$.  Note that rightmost sections starting
at $S^1_\infty(\lambda)$ and leftmost sections starting at
$S^1_\infty(\mu)$ also do not cross.

\subsection{Monotone relations between comparable leaves}

We will give two different points of view on how the leftmost sections of
$C$ give a monotone relation from $S^1_\infty(\mu)$ to
$S^1_\infty(\lambda)$.  

\begin{figure}
\psfrag{p}{$p$}
\begin{center}
  \includegraphics[scale=0.5]{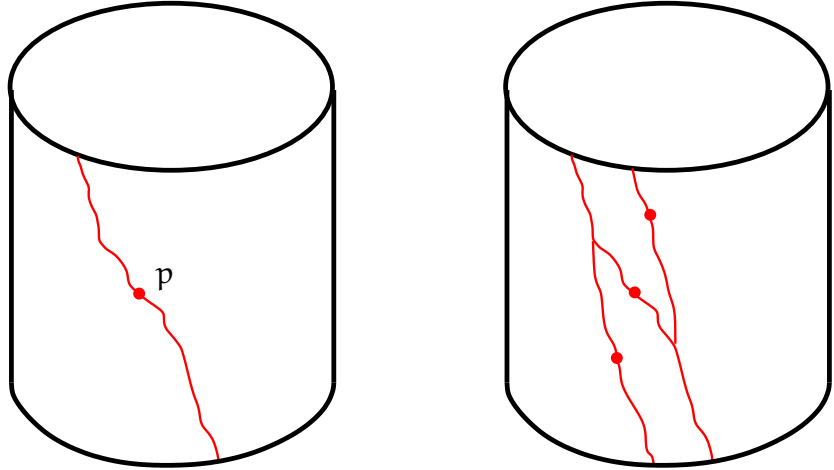}
\end{center}
\caption{Some special sections, showing how they can coalesce.}\label{sections_fig}
\end{figure}

Given a $\nu \in I$ and a $p$ in $S^1_\infty(\nu)$, we define the
\emph{special section} of $C$, denoted $\tau_p$, to be the section
which is the leftmost section going up from $S^1_\infty(\nu)$ to
$S^1_\infty(\lambda)$, and the rightmost section going downward from
$S^1_\infty(\nu)$ to $S^1_\infty(\mu)$ (see
Figure~\ref{sections_fig}).
By the above, two special sections do not cross, though they may
coalesce.  Let $\mathscr{S}$ denote the set of all special sections,
that is, all sections of $C$ which are of the form $\tau_p$ for at
least one $p$.  The paths in $\mathscr{S}$ give us a singular
foliation of $C$ of the type discussed in Lemma~\ref{monotone_param}.
By that lemma, the singular foliation can be split open to give a
monotone relation from $S^1_\infty(\mu)$ to $S^1_\infty(\lambda)$.

Another point of view on this monotone relation is the following,
which is a warm up for the general construction of the universal
circle.  Because $\mathscr{S}$ consists of non-crossing paths, it has
a natural circular order as follows (see the proof of
Theorem~\ref{thm_ord_cataclysms} for the definition of a circular
order).  Consider three distinct sections
\[
(\tau_{p_1}, \tau_{p_2}, \tau_{p_3}).
\]
The non-crossing property guarantees that if $(\tau_{p_1}(\nu),
\tau_{p_2}(\nu), \tau_{p_3}(\nu) )$ are clockwise ordered for some
$\nu \in I$, then for any other $\nu'$, the triple $(\tau_{p_1}(\nu')$,
$\tau_{p_2}(\nu')$, $\tau_{p_3}(\nu') )$ is either clockwise ordered or
degenerate.  Thus, we say that $(\tau_{p_1}, \tau_{p_2}, \tau_{p_3}) $
are clockwise ordered if for some $\nu$ the points $(\tau_{p_1}(\nu)$,
$\tau_{p_2}(\nu)$, $\tau_{p_3}(\nu) )$ are clockwise ordered.  (Actually,
one can have three distinct $\tau_{p_i}$ which intersect each circle
in only two points, but such paths still have a natural circular order
by perturbing the $\tau_{p_i}$ slightly to disjoint paths).

Let $\bar{\mathscr{S}}$ be the closure of $\mathscr{S}$ in the space
of sections $I \to C$.  It is easy to check that every $\tau \in
\bar{\mathscr{S}}$ is admissible, and that any two sections in
$\bar{\mathscr{S}}$ do not cross.  We claim that $\bar{\mathscr{S}}$
is isomorphic, as a circularly ordered set, to a circle.  For this, we
just need to check that $\bar{\mathscr{S}}$ is compact under the
topology induced by the circular order.  But this is the case because
the topology induced by the order is the same as the one induced as a
subspace of the space of all sections.  Henceforth, we will denote the
circle $\bar{\mathscr{S}}$ by $S^1_{\mu \lambda}$.  

For any $\nu \in I$, we have a monotone map from $S^1_{\mu \lambda}$
to $S^1_\infty(\nu)$ by evaluation of sections.  In particular, we
have our promised monotone relation by looking at the restrictions to
$\mu$ and $\lambda$.  The circle $S^1_{\mu \lambda}$ should be thought
of as the universal circle of the foliation restricted to $I$.

\subsection{The universal circle for $\R$-covered foliations}

Before doing the general case, let's construct the universal circle in
the case where the leaf space $L = \R$.  For $p \in S^1_\infty(\nu)$
we can define the \emph{special section}, $\tau_p$, of $E_\infty$ as
above (that is, it's the leftmost section upward from $p$, and the
rightmost section downward).  Let $\mathscr{S}$ denote the set of all
such special sections, which has a natural circular order.  Moreover,
we have surjective monotone maps $\mathscr{S} \to S^1(\nu)$ for each
$\nu \in L$ via evaluation.

We claim that $\mathscr{S}$ is invariant under the action of
$\pi_1(M)$ on $E_\infty$.  The only thing to worry about is that if $\F$ is
not co-orientable, then the action of $\pi_1(M)$ on $\til{\F}$ can
reverse the co-orientation we used to define leftmost.  However, since
$M$ is orientable, if $\gamma \in \pi_1(M)$ flips the co-orientation
then $\gamma$ also reverses the orientation of all the circle fibers
of $E_\infty$.  This has the affect of exchanging leftmost upwards
paths with rightmost downwards paths and vice versa.  Thus
$\mathscr{S}$ is preserved by $\pi_1(M)$.

To construct the universal circle, we just need to complete
$\mathscr{S}$ into a circle.  Up to homeomorphism, there is a unique
way to embed $\mathscr{S}$ into a circle $S^1$ so that the embedding
respects the circular order and is continuous with respect to the
topology on $\mathscr{S}$ induced by its order.  The closure of the
image of $\mathscr{S}$ might omit some gaps, which we can collapse to
get a monotone map from $\mathscr{S}$ to a circle $S^1_\un$, which is
the promised \emph{universal circle}.  Because no $\tau_p$ is isolated
to either side in $\mathscr{S}$, it follows that $\mathscr{S}$
actually embeds in $S^1_\un$.

The action of $\pi_1(M)$ on $\mathscr{S}$ extends naturally to an
action on $S^1_\un$.  Likewise, we have natural monotone maps from
$S^1_\un$ to each $S^1_\infty(\nu)$.  It easy to check that $S^1_\un$
has all the properties of Definition~\ref{univ_circle_def}, and so
we've proved Theorem~\ref{universal_circle_exists} in the case where
$L = \R$.

\subsection{Turning corners}

To build a universal circle in general, we need to figure out how to
relate the circles at infinity of incomparable leaves.  Recall that
two incomparable leaves $\lambda_1$ and $\lambda_2$ are part of the same
cataclysm if they are both limits of a sequence $\nu_i$ of comparable
leaves (see Section~\ref{sec_non_tight}).  The following lemma is the
key to constructing special sections in the non--$\R$-covered case.
\begin{Lemma}[Turning corners]\label{turning_corners}
  Let $\lambda_1$ and $\lambda_2$ be a pair of incomparable leaves of the same
  cataclysm. Then there is a canonical monotone relation between
  $S^1_\infty(\lambda)$ and $S^1_\infty(\mu)$.
\end{Lemma}
\begin{proof}
  Roughly, the monotone relation is the inverse limit of the circles
  at infinity in the approximating leaves of the cataclysm.  We'll
  assume that $\lambda_1$ and $\lambda_2$ are a limit from below.
  That is, there are intervals in the leaf space $I_i = [\mu,
  \lambda_i]$, whose intersection is $I = [\mu, \lambda_1) = [\mu,
  \lambda_2)$ (this is the situation shown inside the circle in
  Figure~\ref{between}).  Let $C_i$ denote the cylinder $E_\infty
  \vert_{I_i}$, and $C$ the half-open cylinder $E_\infty
  \vert_{I}$.  Let $\mathscr{M}_i$ be the set of markers heading
  downward from $S^1_\infty(\lambda_i)$ into $C_i$.  Essentially, the
  needed monotone relation is obtained by completing the set of
  markers $\mathscr{M} = \mathscr{M}_1 \cup \mathscr{M}_2$ into a
  circle.
  
  First we claim that if $m_1 \in \mathscr{M}_1$ and $m_2 \in
  \mathscr{M}_2$ then $m_1$ and $m_2$ are disjoint in $C$.  If
  not, as in the proof of Lemma~\ref{marker_crossing}, we can join
  points $x_i \in \lambda_i$ by a path $\gamma$ of length less than a
  separation constant $\epsilon$.  A leaf $\mu \in I$ close to the
  $\lambda_i$ intersects $\gamma$ in a pair of points $y_i$, one close
  to each $x_i$.  The distance between the $y_i$ in $\mu$ can be made
  arbitrary large by choosing $\mu$ close enough to the $\lambda_i$.
  As the $y_i$ are distance less than $\epsilon$ in $\til{M}$, this
  violates that $\epsilon$ is a separation constant.  

  We can define a map $f_i \maps \mathscr{M} \to
  S^1_\infty(\lambda_i)$ as follows.  For a marker $m \in
  \mathscr{M}$, look at the (possibly half-open) marker $m \cap C_i$.
  Denote by $\overline{m}$ the closure of $m \cap C_i$ in $C_i$.  If
  $m$ is in $\mathscr{M}_i$, then of course $\overline{m}$ is just $m$.
  In general, the density of markers implies that $\overline{m}$ is
  an interval transverse to the circle fibers of $C_i$.  We define
  $f_i(m)$ to be the point of intersection of $\overline{m}$ with
  $S^1_\infty(\lambda_i)$.  Consider pairs of markers in
  $\mathscr{M}_1$ and $\mathscr{M}_2$.  Because the $\lambda_i$ are
  incomparable, the intersection of these pairs of markers with any
  $S^1_\infty(\nu)$ in $C$ are unlinked.  This implies that for $i \neq j$,
  $f_i(\mathscr{M}_j)$ consists of a single point $p_i \in
  S^1_\infty(\lambda_i)$.  The set $\mathscr{M}$ has a circular order
  coming from intersecting markers with $C$.  The
  completion/blow-down of  $\mathscr{M}$ is a circle
  $S^1_{\lambda_1 \lambda_2}$ which is the needed relation.
  
  Explicitly, we can construct $S^1_{\lambda_1 \lambda_2}$ by
  cutting each $S^1_\infty(\lambda_i)$ at $p_i$ and gluing the
  resulting intervals together to form a circle.   The map of
  the relation $f_i \maps S^1_{\lambda_1 \lambda_2} \to
  S^1_\infty(\lambda_i)$ is just the identity on the interval
  coming from  $S^1_\infty(\lambda_i)$, and maps the interval
  coming from $S^1_\infty(\lambda_j)$ to the point $p_i$.
\halmos \end{proof}

\begin{rmk}\label{distinct_continuation}
  The construction in Lemma~\ref{turning_corners} has the following
  important property.  Take two points $q_i \in \lambda_i$ on a pair
  of incomparable leaves situated as above.  Consider the rightmost
  (downward) section $\tau_i$ on $C_i$.  We claim that these sections
  $\tau_i$ differ on $C$.  Pick a sequence of leaves $\mu_j$
  converging to the $\lambda_i$.  Denote by $\mathscr{M}_i(j)$ the
  intersection of $\mathscr{M}_i$ with $\mu_j$.  The subsets
  $\mathscr{M}_i(j)$ are nonempty for large $j$, and they are
  \emph{disjoint}.  Moreover, no pair of points in $\mathscr{M}_1(j)$
  links a pair of points in $\mathscr{M}_2(j)$.  Because the $\tau_i$
  are rightmost sections, it follows that either $\tau_1(j) =
  \tau_1(\mu_j)$ is separated from $\tau_2(j)$ for some $j$ by pairs of
  points in $\mathscr{M}_1(j)$, or else $(m_1(j),\tau_1(j),m_2(j))$ is
  anticlockwise ordered for each pair of markers $m_1 \in
  \mathscr{M}_1,m_2 \in \mathscr{M}_2$.  Conversely, either
  $\tau_2(j)$ is separated from $\tau_1(j)$ for some $j$ by pairs of
  points in $\mathscr{M}_2(j)$, or else $(m_2(j),\tau_2(j),m_1(j))$ is
  anticlockwise ordered for each pair of markers $m_1 \in
  \mathscr{M}_1,m_2 \in \mathscr{M}_2$. It follows that either the
  $\tau_i(j)$ are separated from each other by $\mathscr{M}_i(j)$, or
  else for large $j$, there are $m_1(j),m_2(j)$ such that
  $(m_1(j),\tau_1(j),m_2(j),\tau_2(j))$ are anticlockwise ordered; in
  particular, $\tau_1(j)$ and $\tau_2(j)$ are distinct in this case
  too.  Note that the $\tau_i(j)$ only fail to be separated from each
  other by pairs of points in some $\mathscr{M}_i(j)$ if and only if
  the points $q_i$ are precisely the special points $p_i$ used in the
  explicit construction of $S^1_{\lambda_1 \lambda_2}$.

\end{rmk}

\subsection{Special sections in general} 
\begin{figure}[hbt]
  \psfrag{a}{$\lambda$}
  \psfrag{b}{$\mu$}
  \begin{center}
    \includegraphics[scale=0.5]{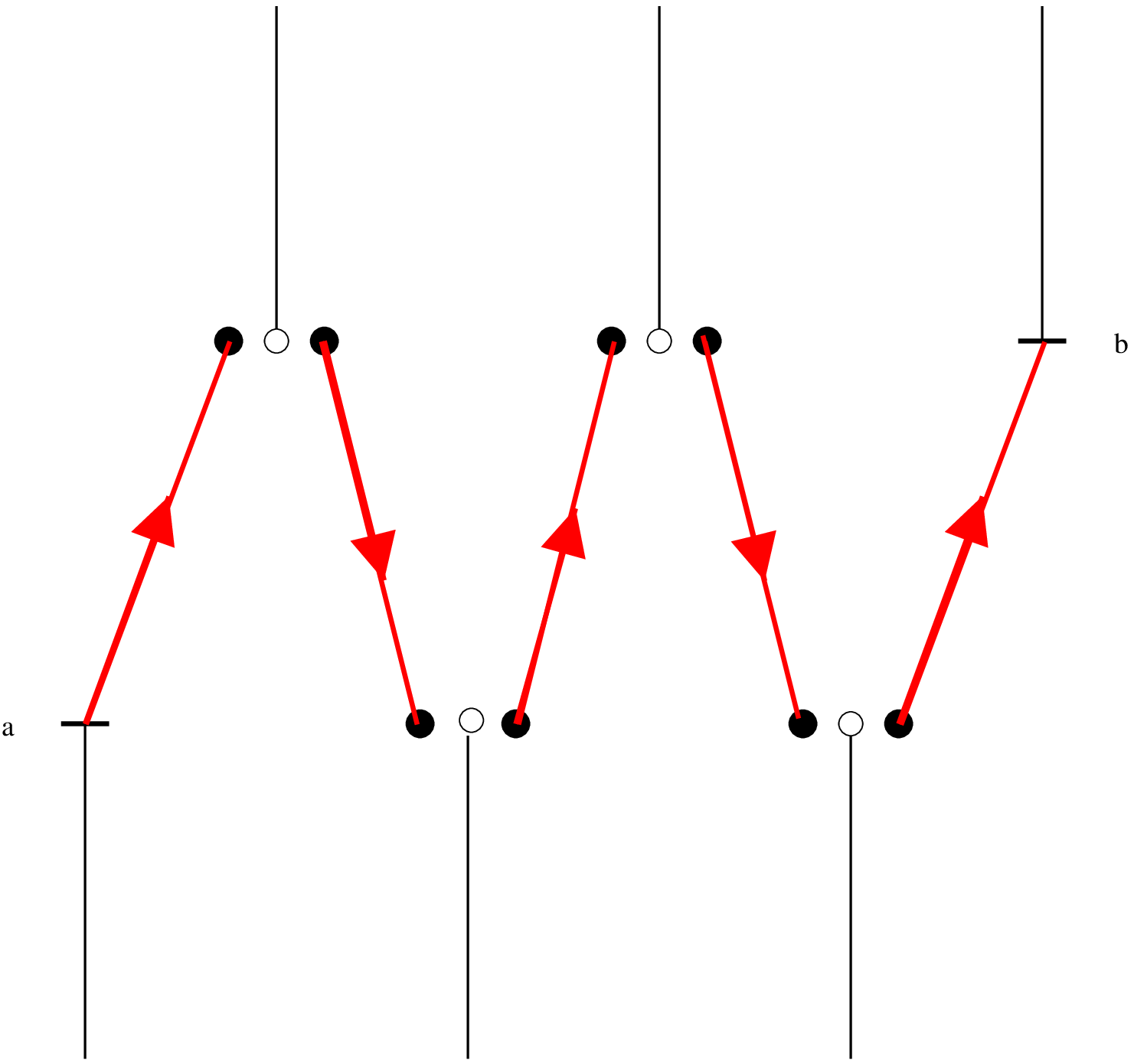}
  \end{center}
\caption{The path $\gamma_{\lambda \mu}$.}\label{section}
\end{figure}

For a point $p \in S^1_\infty(\lambda)$, we will define the associated
\emph{special section}, $\tau_p$, as follows.  For a leaf $\mu$ in
$\til{\F}$ there is a unique minimal broken ``path'' that joins
$\lambda$ to $\mu$ in $L$ consisting of finitely many ordered
intervals and turns.  That is, there is a sequence
\[
\lambda = \mu_0,\mu_1,\mu_2,\mu_3,\dots,\mu_n = \mu
\]
where $(\mu_{2i}, \mu_{2i+1})$ are comparable, and $(\mu_{2i+1},
\mu_{2i+2})$ are limit leaves of the same cataclysm (see
Figure~\ref{section}).
Call this path $\gamma_{\lambda \mu}$.  We define $\tau_p$ on the
restriction of $E_\infty$ to $\gamma_{\lambda \mu}$ by starting at $p$
in $S^1_\infty(\lambda)$ and then extending out along $\gamma_{\lambda
  \mu}$ in the following way.  On the segments of $L$ joining
$\mu_{2i}$ to $\mu_{2i+1}$, we use the leftmost section when
$\gamma_{\lambda \mu}$ is heading upwards and the rightmost section
when $\gamma_{\lambda \mu}$ is heading downwards.  For the turns
between $\mu_{2i+1}$ and $\mu_{2i+2}$, we use the monotone relation of
Lemma~\ref{turning_corners}, which associates a point in
$S^1_\infty(\mu_{2i +2})$ with $ S^1_\infty(\mu_{2i + 1})$, to continue
$\tau(p)$.  This then defines the value of $\tau_p$ at $\mu$.

In order to construct the universal circle we need to understand how
$\tau_p$ depends on $p \in S^1_\infty(\lambda)$.  If $q$ is another
point in the same $S^1_\infty(\lambda)$, then $\tau_p = \tau_q$ for
any leaf $\mu$ which is incomparable to $\lambda$.  Moreover, for $q$
in an arbitrary $S^1_\infty(\mu)$, the two sections $\tau_p$ and
$\tau_q$ can differ on only part of $L$ in the following sense.  Let
$\bar{\gamma}_{\lambda \mu}$ denote the union of $\gamma_{\lambda
  \mu}$ with all $\nu$ which are part of the same cataclysm as some
$(\mu_{2i + 1}, \mu_{2i + 2})$. We say that a line $A \subset L$ is
between $\lambda$ and $\mu$ if it intersects $\bar{\gamma}_{\lambda
  \mu}$.  Then:
\begin{Lemma}\label{diffing_point}
  Suppose $\tau_p$ and $\tau_q$ are sections which differ on a leaf
  $\nu \in L$.  Then $\nu$ lies on a line $A$ which lies
  between $\lambda$ and $\mu$.
\end{Lemma}
\begin{proof}
  Suppose that $\nu$ does not lie on a line between $\lambda$ and $\mu$.
  Then the paths $\gamma_{\lambda \nu}$ and $\gamma_{\mu \nu}$
  intersect in a nontrivial subpath $\gamma$ ending at $\nu$.  Note
  that the orientations of $\gamma_{\lambda \nu}$ and $\gamma_{\mu
    \nu}$ agree on $\gamma$, and that $\gamma$ contains at least one
  turn.  The sections $\tau_p$ and $\tau_q$ coalesce as they pass
  through that turn, and so $\tau_p(\nu) = \tau_q(\nu)$.  
\halmos \end{proof}
On the other hand, we have:
\begin{Lemma}\label{line_differ}
  Let $\tau_p$ and $\tau_q$ be two distinct sections, and $A$ a line in
  $L$ which lies between $\tau_p$.  Then $\tau_p$ and $\tau_q$ differ on
  $A$.
\end{Lemma}
\begin{proof}
  If $\lambda$ and $\mu$ are comparable, then as $\tau_p$ and $\tau_q$
  are distinct they must differ on the interval $I = [\lambda, \mu]$
  in $L$.  As the line $A$ contains $I$, the two sections differ on
  $A$.  So now assume that $\lambda$ and $\mu$ are incomparable.
  \begin{figure}[htb]
    \psfrag{a}{$\lambda$}
    \psfrag{b}{$\mu$}
    \psfrag{c}{$I$}
    \psfrag{d}{$A$}
    \begin{center}
      \includegraphics[scale=0.3]{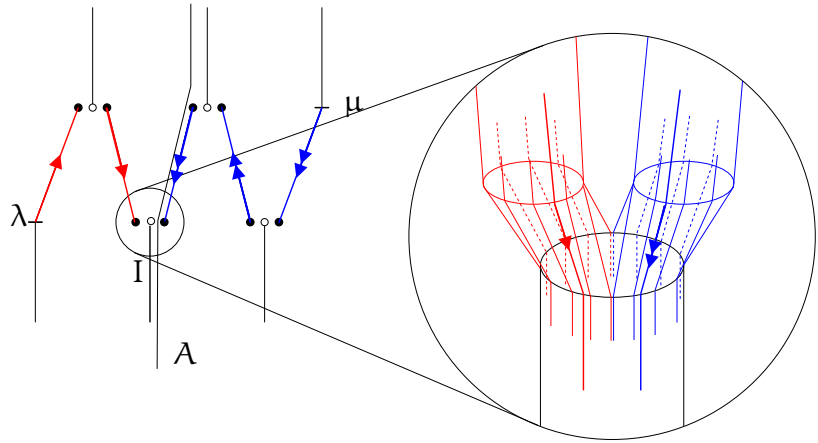}
      \end{center}
    \caption{Where the sections can differ.} \label{between}
  \end{figure}
  As in Figure~\ref{between}, $A$ must pass through a turn of
  $\gamma_{\lambda \mu}$.  Thus we have an interval $I$ of $A$ where
  the two sections come in via different branches of the cataclysm.
  So the situation looks like the magnified picture in
  Figure~\ref{between}.  As noted in
  Remark~\ref{distinct_continuation}, the way the monotone relation
  works for incomparable leaves implies that the two sections differ
  on $I$.  
\halmos \end{proof}

Now let $\mathscr{S}$ denote the set of all special sections $\tau_p$.
We want to put a natural circular order on $\mathscr{S}$ by declaring
that a triple $(\tau_{p_1}, \tau_{p_2}, \tau_{p_3})$ is clockwise
ordered if the restrictions of the $\tau_{p_i}$ to some line $A \subset
L$ are clockwise ordered.  To know  that this makes sense, we need to
check:
\begin{Lemma} 
  Given three distinct sections $(\tau_{p_1}, \tau_{p_2},
  \tau_{p_3})$, there exists a line $A$ so that the restriction of the
  $\tau_i$ to $A$ has a non-degenerate circular order.  Moreover, this
  circular order is independent of the choice of $A$.
\end{Lemma}
\begin{proof}
  Let $\lambda_i$ be the leaf where $p_i \in S^1_\infty(\lambda_i)$.
  The \emph{midpoint} of the $\lambda_i$ is defined analogously to the
  midpoint of three points in an $\R$-tree.  More precisely, we let
  $\Gamma$ be the union of the three paths $\gamma_{\lambda_i
    \lambda_j}$.  The midpoint is constructed by making $\Gamma$
  Hausdorff by amalgamating the cataclysms, taking the midpoint in the
  resulting tree, and then pulling back that midpoint to $L$.  The
  midpoint consists of either a single point in $L$ or several points
  of the same cataclysm.  Given Lemma~\ref{diffing_point}, if $\nu$ is
  a leaf where all the $\tau_{p_i}$ differ, then for each pair $(i,j)$
  there is a line $A_{ij}$ containing $\nu$ which intersects
  $\bar{\gamma}_{\gamma_i \gamma_j}$.  So a natural place to look for
  the needed line is near the midpoint of the $\lambda_i$.
  \begin{figure}
    \psfrag{A}{(a)}
    \psfrag{B}{(b)}
    \psfrag{C}{(c)}
    \psfrag{A}{(a)}
    \psfrag{11}{$\lambda_1$}
    \psfrag{21}{$\lambda_2$}
    \psfrag{31}{$\lambda_3$}
    \psfrag{a2}{$\lambda_1$}
    \psfrag{22}{$\lambda_2$}
    \psfrag{32}{$\lambda_3$}
    \psfrag{a1}{$I$}
    \psfrag{L}{$A$}
    \begin{center}
      \includegraphics[scale=0.5]{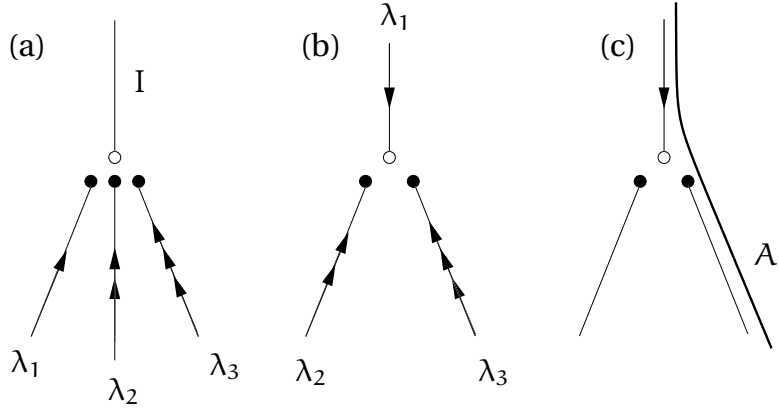}
    \end{center}
    \caption{Some possible cases.}\label{midpoint}
  \end{figure}
  There are several cases.  First, consider the case where the
  midpoint is several points of the same cataclysm.  There are two
  subcases corresponding to Figure~\ref{midpoint}(a) and (b).  In case
  (a), the proof of Lemma~\ref{line_differ} shows that the three
  sections differ on the interval $I$ shown.  In case (b), first
  notice that $\tau_{p_2}$ and $\tau_{p_3}$ differ on the interval
  $I$.  Therefore, $\tau_{p_1}$ differs from at least one of
  $\tau_{p_2}$ and $\tau_{p_3}$ on $I$.  Without loss of generality,
  assume that $\tau_{p_1} \neq \tau_{p_2}$.  Take $A$ to be the line
  shown in Figure~\ref{midpoint}(c).  As $A$ is between $\lambda_1$
  and $\lambda_3$, the sections $\tau_{p_1}$ and $\tau_{p_3}$ differ
  on $A$.  As $I \subset A$, we know that the $\tau_{p_i}$ are
  distinct on $A$, and so have a  non-degenerate circular order.
  
  In the case that the midpoint is a single point $\mu$, take any line
  $A$ containing $\mu$.  The midpoint $\mu$ lies in every
  $\gamma_{\lambda_i \lambda_j}$, and so $A$ lies between each pair
  $(\lambda_i, \lambda_j)$.  By Lemma~\ref{line_differ}, the $\tau_i$
  restricted to $A$ are distinct, and so have a non-degenerate
  circular order.
  
  It remains to prove that the circular order is independent of the
  choice of the line $A$.  Because of Lemma~\ref{diffing_point}, this
  is simple to check in each of the above cases, and we leave this to
  the reader.
\halmos \end{proof}

Now that we've defined special sections in general and shown they have
a circular order, the proof of Theorem~\ref{universal_circle_exists}
follows as in the $\R$-covered case by completing $\mathscr{S}$ into a
circle.  The only difference is that now we need to check Property
(\ref{assoc_gaps}) of Definition~\ref{univ_circle_def}, but this is
clear from the construction in Lemma~\ref{turning_corners}, and
from Remark~\ref{distinct_continuation} immediately following it.

\begin{rmk}
  Note that there is no claim of {\em uniqueness} for the universal
  circles constructed by this or other methods. Even if one asks for
  minimal universal circles, it is not {\it a priori} clear that they
  should be unique. On the other hand, for $\R$-covered foliations or
  those with one-sided branching, there is a unique minimal universal
  circle \cite{Calegari2000,Calegari:onesided}.
\end{rmk}

\subsection{Faithfulness of the action}
Finally, we conclude this section by proving
Theorem~\ref{universal_circle_action_faithful}, namely that the action
of $\pi_1(M)$ on $S^1_\un$ is faithful provided $M$ is atoroidal.  If
$\F$ is $\R$-covered, then since $M$ is atoroidal there is a
transverse pseudo-Anosov flow which can be split open to a pair of
very full genuine laminations $\Lambda^\pm$ transverse to each other
and to $\F$. These laminations are covered by genuine laminations
$\til{\Lambda}^\pm$ which arise from a pair of invariant laminations
$\Lambda^\pm_\un$ of $S^1_\un$.  If $\gamma$ acts as the identity on
$S^1_\un$, it acts as the identity on the leaf space of either genuine
lamination, and therefore it fixes the singular flowlines of the
pseudo-Anosov flow (see \cite{Calegari2000} Corollary 5.3.16); in
particular, it acts as a nontrivial translation along these flowlines.
But the dynamics of the flow expands the transverse measure on one
singular lamination and contracts the transverse measure on the other,
and therefore the action on the leaf spaces of the singular foliations
is nontrivial, giving a contradiction.  So we can assume that the leaf
space $L$ has branching (in this case, we will not need the assumption
that $M$ is atoroidal).

Let $K$ be the kernel of the map $\pi_1(M) \to \homeo(S^1_\un)$.  Let
$\lambda$ be any leaf of $\til{\F}$, and $k$ a non-identity element in
$K$.  We claim that $k \lambda$ and $ \lambda$ are comparable and
distinct.  Note that $k$ can't fix $\lambda$ setwise, because if it
did it would act identically on $S^1_\infty(\lambda)$ and hence
identically on $\lambda$.  If $\lambda$ and $k \lambda$ are
incomparable, we know that the maps $S^1_\un \to S^1_\infty(\lambda)$
and $S^1_\un \to S^1_\infty( k \lambda)$ blow down different gaps.  As
$k$ acts identically on $S^1_\un$, this is impossible.

Fix a branch point $\lambda$ in $L$.  By our above observations, the
set $K \lambda$ is infinite and contained in a line $A \subset L$.
Now consider a leaf $\mu$ which is in the same cataclysm as $\lambda$.
Again, $K \mu$ is contained in a line $B$.  For each $k$, the leaf $k
\mu$ is in the same cataclysm as $k \lambda$.  So there are infinitely
many pairs of non-separable points $(k \lambda, k \mu)$ where one
point is in $A$ and other in $B$.  This is not possible for two lines
in $L$, and so we have a contradiction.  Thus the action on $S^1_\un$
is faithful in the branching case.  This completes the proof of
Theorem~\ref{universal_circle_action_faithful}.

\section{Group actions on $1$-manifolds and left orderings}\label{sec_orderings}

In this section, we discuss the algebraic interdependence between the
existence of actions of a group $G$ on various $1$-manifolds.  The
main reason for this is to improve the usefulness of the
(non)existence of such a group action as a criterion for the existence
of foliations and laminations.  In particular, we want to promote an
action of $\pi_1(M)$ on $S^1$ to an action on $\R$.  This is because
it is much easier to decide whether a group $G$ admits a a faithful
action on $\R$ than on $S^1$.

We'll begin with a more algebraic criterion for a group $G$ to act
faithfully on $\R$.  A group $G$ is \emph{left-orderable} if there is
a total order on $G$ which is invariant under left multiplication. If
$G$ has a faithful representation into $\homeo^+(\R)$, then $G$ is
left-orderable for the following reason.  Pick an infinite sequence of
points $p_1,p_2,\dots \in \R$ such that the intersection of the
stabilizers of these points in $G$ is trivial.  Then say $\alpha >
\beta$ if for the smallest $i$ where one of $\alpha$ and $\beta$ does
not fix $p_i$, we have $\alpha(p_i) > \beta(p_i)$. Conversely, for a
finitely generated group $G$, it is not hard to show that the
existence of a left-invariant order gives rise to a faithful
representation in $\homeo^+(\R)$, so these conditions are actually
equivalent (see e.g.~\cite{Ghys:circle}).  For actions on circles,
one can use the analogous notion of a circular order to the same
effect.

\subsection{Lifting representations}\label{lifting}

Given a group $G$ and a representation 
\[
\rho \maps G \to \homeo(S^1),
\]
when does $\rho$ lift to a representation $\til{\rho} \maps G \to
\homeo(\R)$? Here we want the lift to respect the universal covering
map $\pi \maps \R \to S^1$ in the sense that
\[
\pi \circ \til{\rho}(g) = \rho \circ \pi(g) \mtext{for all $g$ in $G$.}
\]

The topological group $\homeo(S^1)$ has two connected components; the
component of the identity is the subgroup $\homeo^+(S^1)$ of
orientation preserving homeomorphisms.  If $H^1(G;\Z/2\Z)=0$, every
homomorphism from $G$ to $\homeo(S^1)$ has image in $\homeo^+(S^1)$.

There is a universal central extension of $\homeo^+(S^1)$
\[
1 \to \Z \to \til{\homeo^+}(S^1) \to \homeo^+(S^1) \to 1
\]
where the group $\til{\homeo}(S^1)$ is the group of  periodic
  homeomorphisms of $\R$ with period $1$; i.e.~those that commute
with the translation $Z \maps t \to t+1$.

Since $\til{\homeo^+}(S^1)$ is the universal central extension of
$\homeo^+(S^1)$, the obstruction to lifting a homomorphism $\rho \maps
G \to \homeo^+(S^1)$ is an element of $H^2(G;\Z)$ called the \emph{
  Euler class} $e(\rho)$ of $\rho$ (for more, see \cite{Ghys:circle}).
We will give a geometric interpretation of $e(\rho)$ later, but for
the moment we study the implications for $G = \pi_1(M)$ where $M$ is a
rational homology sphere.

\begin{Theorem}\label{dihedral_cover}
  Let $M$ be a rational homology sphere and $\rho \maps \pi_1(M) \to
  \homeo^+(S^1)$ a faithful homomorphism.  Then the commutator
  subgroup $[\pi_1(M), \pi_1(M)]$ is left-orderable.  
\end{Theorem}
\begin{proof}
  Let $G = \pi_1(M)$.  By assumption $H_1(M;\Z)$ is finite, so that
  $e(\rho) \in H^2(G;\Z) \iso H_1(M; \Z)$ is torsion.  We will show
  that $\rho$ restricted to $[G, G]$ lifts to $\til{\homeo^+}(S^1)$.

For each $g \in
  G$, there is a lift of $\rho(g)$ to $\til{\rho}(g) \in
  \til{\homeo^+}(S^1)$.  This lift is not unique, and distinct lifts
  differ by powers of the generator $Z$ of the center of
  $\til{\homeo^+}(S^1)$. If $c \in [G,G]$, we can express $c$ as a
  product of commutators
  \[
  c = \prod_i [g_i,h_i].
  \]
  Consider the element in  $\til{\homeo^+}(S^1)$ given by
  \[
  \til{\rho}(c) = \prod_i [\til{\rho}(g_i),\til{\rho}(h_i)].
  \]
  Since $Z$ is central, this depends only the description of $c$ as a
  product of commutators, not on the choices of the
  $\{\til{\rho}(g_i), \til{\rho}(h_i)\}$.  We will show that
  $\til{\rho}(c)$ is independent of the choice of expression as a
  product of commutators.  Then $\til{\rho} \maps [G, G] \to
  \til{\homeo^+}(S^1)$ will be a homomorphism which is a lift of
  $\rho$ as required.

  If we have another description of $c$ as a product of commutators
  \[
  c = \prod_i [g_i',h_i']
  \]
  then the relation
  \[
  \text{id} = \prod_i [g_i,h_i]\left(\prod_j [g_i',h_i']\right)^{-1}
  \]
  determines a map from a surface $\Sigma \to M'$, and therefore an element 
  $[\Sigma] \in H_2(G;\Z)$. To see that $\til{\rho}$ is
  well-defined on $[G,G]$ it suffices to show that
  \[
  Z^n = \prod_i [\til{\rho}(g_i),\til{\rho}(h_i)] \left(
    \prod_j [\til{\rho}(g_j'),\til{\rho}(h_j')] \right)^{-1}
  \]
  is trivial in $\til{\homeo^+}(S^1)$. But by the definition of the Euler class,
  we have that
  \[
  n = e([\Sigma])
  \]
  is zero, because $e$ is torsion. So the representation lifts to
  $\til{\homeo^+}(S^1)$ on $[G,G]$, and we are done.
\halmos \end{proof}

In case of a general action on $S^1$ we have:
\begin{Theorem}\label{leftorder_subgroup_gen}
  Let $M$ be an irreducible rational homology sphere, and 
  \[ 
  \rho \maps \pi_1(M) \to \homeo(S^1)
  \] 
  a faithful homomorphism.  Let $K \subset
  \pi_1(M)$ be the kernel of the homomorphism from $\pi_1(M)$ to
  $\Z/2\Z$ defined by the orientation of $\rho$.  Then the subgroup
  $[K, K]$ is left-orderable.
\end{Theorem}
\begin{proof}
  Consider the manifold $M'$ corresponding to $K$.  If $M'$ is also a
  rational homology sphere, then we are done by the preceding theorem.
  If instead $H^1(M) \neq 0$, then $K = \pi_1(M')$ itself is
  left-orderable for reasons that have nothing to do with the
  existence of $\rho$, namely the following theorem of from \cite[Cor.
  3.4]{BoyerRolfsenWiest2002}:
  \begin{Theorem}[Boyer, Rolfsen, Wiest]\label{hom_to_order}
    A compact, orientable, irreducible 3-man\-i\-fold with $H^1(M) \neq 0$
    has left-orderable fundamental group.
  \end{Theorem}
  In either case, we're done.
\halmos \end{proof}

\begin{rmk}
  The necessity of passing to a dihedral cover in general in
  Theorem~\ref{leftorder_subgroup_gen} reflects the fact that finite
  dihedral groups act faithfully on $S^1$.  For instance, let $M$ be a
  $3$-manifold with a dihedral cover $N \to M$, where
  $\pi_1(M)/\pi_1(N) = D_n$, such that $\pi_1(N)$ is left-orderable.
  Let $\rho \maps \pi_1(M) \to \homeo(S^1)$ have image the standard
  action of the dihedral group $D_n$ on $S^1$. The kernel is exactly
  $\pi_1(N)$, and we can modify $\rho$ to a monotonely equivalent {\em
    faithful} action $\rho'$ by blowing up the (finite) orbit of some
  point in $S^1$, inserting a faithful action of $\pi_1(N)$ at some
  blown up interval, and transporting this action around the orbit by
  $D_n$.
\end{rmk}

\begin{Corollary}
  Let $M$ be an orientable atoroidal $3$-manifold which contains a
  taut foliation, a tight essential lamination with solid torus guts,
  or a pseudo-Anosov flow. Then there is a finite cover of $M$ with
  left-orderable fundamental group.  The covering group is either an
  abelian group, or a $\Z/2\Z$ extension of an abelian group.
\end{Corollary}
\begin{proof}
  By the results of previous sections, in every case of the hypothesis
  of the theorem, there is a faithful representation of $\pi_1(M)$
  into $\homeo(S^1)$.  If $H^1(M;\Z) \ne 0$, then $\pi_1(M)$ is
  left-orderable by Theorem~\ref{hom_to_order}.  Otherwise, the result
  follows from Theorem~\ref{leftorder_subgroup_gen}.
\halmos \end{proof}

\subsection{The Euler class and multisections}\label{multisections}

A more geometric way to think about Euler classes is via foliated
bundles and multisections.  A representation $\rho \maps \pi_1(M) \to
\homeo^+(S^1)$ determines a foliated $S^1$ bundle $E$ over $M$, by
taking the quotient of the trivial bundle $\til{M} \times S^1$ by
$\pi_1(M)$ acting via
\[
\gamma(m,\theta) = (\gamma(m),\rho(\gamma)(\theta)).
\]
Another interpretation of the Euler class $e$ of $\rho$ is the
obstruction to finding a section of this foliated bundle. For any
integer $n$, the class $ne$ is the obstruction to finding an {\em
  order $n$ multisection} of this bundle.  The existence of such a
multisection gives a representation of $\pi_1(M)$ on the group of
homeomorphisms of $n$ copies of $\R$, as follows. An order $n$
multisection can be parameterized on a small open set $B \subset M$ by
$n$ sections $s_i \maps B \to E$. The ``labels'' on these sections can
be analytically continued along a path in $M$, but after traveling
around a loop $\gamma \subset M$, the labels are permuted by some
element in the symmetric group on $n$ letters:
\[
\sigma(\gamma) \in \mathsf{S}_n.
\]
Pick a point $p$ in $B$, and consider the cover of the fiber $S^1_p$
by $n$ copies of $\R$,  where the basepoint of the $i$th copy of $\R$
maps to $s_i(p)$.  As we wind around a loop $\gamma$ based at $p$, we
can see how much holonomy transport of some leaf twists in $S^1_p$ {\em
  relative} to the section $s_i$.  But after traveling around the loop,
the section $s_i$ gets relabelled as $s_{\sigma(\gamma)(i)}$.  This
holonomy transport gives an identification of the $i$th copy of $\R$
with the $\sigma(\gamma)(i)$th copy of $\R$, and one gets a
representation of $\pi_1(M)$ in the group $\homeo^+(\R \times \lbrace
1,\dots,n\rbrace)$.

\subsection{Flips and canceling Euler classes}

Let $M$ be a $3$-manifold, and suppose that $M$ contains a taut
foliation, a tight essential lamination with solid torus guts, or a
pseudo-Anosov flow. Then there is a faithful representation $\pi_1(M)
\to \homeo(S^1)$. We would like to find a subgroup $K$ of
$\pi_1(M)$, of as large index as possible, where the representation
lifts to $\homeo^+(\R)$. From what we have said above, it might seem
that the worst case is when $H^1(M;\Z/2\Z) \ne 0$. But there are some
instances where the nature of these representations helps out, using
the flips discussed in Section~\ref{flips_one}.

Suppose $\Lambda$ is a tight essential lamination of $M$ with solid
torus guts, and suppose further that there is a homomorphism $\sigma
\maps \pi_1(M) \to \Z/2\Z$ such that the core of each gut region maps
to the identity. Let $\widehat{M}$ be the $2$-fold cover corresponding
to the kernel of $\sigma$.  Then $\Lambda$ lifts to a tight essential
lamination of $\widehat{M}$ with two solid torus gut complementary
regions for each solid torus gut region of $\Lambda$. By the Filling
Lemma~\ref{filling_lemma}, we can assume the complementary regions to
$\Lambda$ were actually ideal polygon bundles over $S^1$, and the same
is therefore true for $\widehat{\Lambda}$. The complementary regions
to $\widehat{\Lambda}$ come in pairs $C_1^i,C_2^i$ corresponding to
the complementary regions $C^i$ to $\Lambda$. The representation $\rho
\maps \pi_1(M) \to \homeo(S^1)$ restricts to a representation $\rho
\maps \pi_1(\widehat{M}) \to \homeo(S^1)$. If, for each complementary
region $C^i$, we perform a flip of {\em exactly one} of the $C^i_j$,
we get a {\em new} representation $\rho'$ which is evidently both
orientation-preserving and has Euler class $0$. In particular, $\rho'$
lifts to a faithful representation in $\homeo^+(\R)$, and
$\pi_1(\widehat{M})$ is left-orderable.  Unfortunately, we see no
means of exploiting this trick for algorithmic purposes, since the
condition that $\sigma$ be trivial on the core of every gut region
seems hard to verify in practice.

\subsection{Actions on leaf spaces of taut foliations}
For completeness, we include a proof of the following fact mentioned
in the introduction:
\begin{Theorem}\label{fol_action_on_leafspace}
  Let $\F$ be a taut foliation of $M$. Then $\pi_1(M)$ admits a
  faithful action on $L$ without a global fixed point.
\end{Theorem}
\begin{proof}
  Let $K$ be the kernel of the holonomy representation. Then $K$ fixes
  every leaf $\lambda$ of $\til{\F}$. In particular, $K$ is a surface
  group, and is therefore left-orderable, and acts faithfully on $\R$. We
  can insert this action at some point $\lambda \in L$ and transport
  it around by the action of $\pi_1(M)/K$.  Geometrically, we can
  realize such an action by changing the foliation $\F$ by a monotone
  equivalence.
\halmos \end{proof}

\begin{Corollary}
  If $M$ admits a co-orientable $\R$-covered foliation $\F$, then
  $\pi_1(M)$ is left-orderable.
\end{Corollary}

\section{Algorithmic issues}\label{algorithmic_issues}

Let $G$ be a finitely presented group which has a short-lex automatic
structure.  In this section, we describe computational techniques for
proving that $G$ is not left-orderable.  The existence of a short-lex
automatic structure (see \cite{WordProcessingInGroups} for
definitions) implies that the word problem for $G$ is efficiently
solvable.  In fact, there is a set of generators $S$ of $G$ for which
there is a fast algorithm for reducing a word $w$ in $S$ to the
canonical word $w'$ which is lexicographically first among all
shortest words equal to $w$ in $G$.

If $G$ is orderable, consider the positive cone $P = \setdef{g \in
  G}{g > 1}$.  Then $P \cdot P \subset P$ and $G$ is the disjoint
union $P \cup P^{-1} \cup \{1\}$.  Conversely, any $P$ with these two
properties gives rise to a left invariant order via $a > b$ if and
only if $b^{-1} a \in P$.  Let $B(r)$ be the ball in $G$ of radius $r$
about $1$, that is, all words in $S$ of length at most $r$.  For fixed
$r$, we can consider the following
\begin{Question}\label{ballQ}
  Does there exist a $P \subset B(r)$ such that $(P \cdot P) \cap B(r)
  \subset P$ and $B(r)$ is the disjoint union $P \cup P^{-1} \cup
  {1}$?
\end{Question}
If $G$ is orderable, the answer to (\ref{ballQ}) is yes.  If the
answer to (\ref{ballQ}) is no, then $G$ is not orderable.  It is not
hard to show that if $G$ is non-orderable then the answer to
(\ref{ballQ}) is no for large enough $r$.  The idea is that if the answer to
(\ref{ballQ}) is always yes, then one can construct an global positive
cone by taking an inverse limit of the partial positive cones $P(r)
\subset B(r)$, picking at each stage a $P(r) \subset B(r)$ which has
infinitely many extensions to larger $B(R)$.

For any particular $r$, the automatic structure on $G$ combined with
the fact that $B(r)$ is finite means that (\ref{ballQ}) is
algorithmically decidable.  We'll now give a simple recursive
algorithm to do this.  The algorithm follows a similar line to the
proof of Theorem~\ref{weeks_nonorder}, and the reader is advised to read
that proof before proceeding.  Fix $r$ and set $B = B(r)$. The
following recursive procedure, \constructP , has the property that
\constructP(\{ \}) returns \emph{true} if and only if the answer to
(\ref{ballQ}) is yes.
\begin{tabbing}
\quad \= \quad \= \quad \= \kill
\> \constructP ( $P$ such that  $P \subset B$): \\
 \>      \> \textbf{while} $(P \cdot P)  \cap B  \not\subset P$: \\
 \>       \> \>       $P := (P \cdot P) \cap B$ \\
\>\>   \textbf{if} $1 \in P$: \\
\>\>\>        \textbf{return} \textit{false} \\
\> \>     \textbf{if} $B = P \cup P^{-1} \cup \{1\}$: \\
\>\>\>           \textbf{ return} \textit{true} \\
\> \>    $g := $ a shortest word in $B - (P \cup P^{-1} \cup \{1\})$\\
\> \>     \textbf{return} \constructP ( $P \cup \{g\}$) \textbf{or} \constructP ( $P \cup \{g^{-1}\}$ )
\end{tabbing}
Let's get a rough grip on how long this algorithm takes in practice.
A bad case is when $P$ exists as then we have to construct it.  Note
that even if we're handed $P$ by an oracle, it still takes about
$(\#P)^2 = (1/4)(\#B)^2$ multiplications to \emph{check} that $P$
satisfies the conditions in (\ref{ballQ}).  If we were to compute a
multiplication table for $B(r)$ in advance (using $(1/2)(\#B)^2$
multiplications), we could do all further multiplications by table
lookup at essentially no cost.  So in cases where we end up
constructing $P$, we can make the running time $C (\#B)^2$ and this is
roughly the best possible.

However, if the answer to (\ref{ballQ}) is no, in practice one needs
far fewer multiplications and only a few ($\leq 4$) levels of
recursion for the algorithm to finish.  So in practice, a good strategy
seems to be to stop the algorithm if the recursion depth is greater
than $5$, and just assume the final answer would have been yes.

The groups we're interested in are the fundamental groups of
hyperbolic 3-man\-i\-folds.  These do have short-lex automatic
structures.  The problem is that $B(r)$ has exponential growth with
$\#B(r) \approx A C^r$.  For the 2-generator groups we looked at in
Section~\ref{census_exs}, $C$ is usually a little less that $3$.  In
practice, the size of $B(r)$ makes it very difficult to decide
(\ref{ballQ}) if $r$ is bigger than, say, $7$.  Given this, it is
remarkable that among the small volume closed hyperbolic 3-manifolds
in the Hodgson-Weeks census there are quite a few whose fundamental
groups can be shown to be non-orderable by this method (see
Section~\ref{census_exs}).

The algorithm can be easily modified to keep track of the origin of an
element of $P$ as a product of the elements $g$ that have been added
to $P$.  This lets one generate a proof that the group is
non-orderable.  In fact, this was how we found the proof of
Theorem~\ref{weeks_circle}.

\section{The Weeks manifold}

The Weeks manifold is the smallest known hyperbolic 3-man\-i\-fold.
In this section, we will show that its fundamental group cannot act
faithfully on $\R$ or $S^1$.  The Weeks manifold $W$ is the $(5/2,
5/1)$ Dehn surgery on the Whitehead link (our convention is
that $+1$ surgery on a component of this link yields the trefoil).  It
is arithmetic, and its volume is about $0.942707362776$.  Its
fundamental group is
\[
G = \spandef{a, b}{bababAb^2A, ababaBa^2B}
\]
where $A = a^{-1}$ and $B = b^{-1}$.

We'll begin with the easy case of actions on $\R$.  
\begin{Theorem} \label{weeks_nonorder}
   The fundamental group $G$ of the Weeks manifold is non-orderable.
\end{Theorem}

\begin{proof}
  Suppose that $G$ is orderable.  Consider the positive cone 
  \[
  P = \setdef{g\in G}{g > 1}.
   \]
  Then $P \cdot P \subset P$ and $G$ is the
  disjoint union $P \cup P^{-1} \cup \{1\}$.  Switching $>$ around if
  necessary, we can assume that $a > 1$, that is $a \in P$.  We now
  consider the various possibilities.  
    \begin{list}{}{
    \setlength{\labelwidth}{2.2cm}
    \setlength{\leftmargin}{2.5cm}
    \setlength{\labelsep}{0.3cm}
    }
    
  \item[\textbf{Case $b \in P$.}\hfil]  As $P \cdot P \subset P$, we have $ab, bab,
  abab$, etc.~contained in $P$. 
  \begin{list}{}{
      \setlength{\labelwidth}{3cm}
      \setlength{\leftmargin}{3.2cm}
      \setlength{\labelsep}{0.3cm}
      \setlength{\itemindent}{0cm}
  }
  \item[\textbf{Subcase $aB \in P$.}] Then $(abab)\cdot(aB)\cdot
  a\cdot(aB)$ is in $P$.  But $ababaB a^2B = 1$ in $G$, and so $1
  \in P$, a contradiction.
  
  \item[\textbf{Subcase $bA \in P$.}]  Then $(baba)\cdot(bA)\cdot
  b\cdot(bA)$ is in $P$.  But $bababAb^2A = 1$ in $G$, and so $1 \in
  P$, a contradiction.
    \end{list}
 
  \item[ \textbf{Case $B \in P$.}]  Using the relations $R_1$ and $R_2$
    we have 
    \[ 
    BaB^2a^2Ba^2B = b^{-1}R_1^{-1}b \cdot R_2 = 1
    \]
    in $G$.
    But then $BaB^2a^2Ba^2B \in P$, a contradiction.
\end{list}
This shows that such a $P$ does not exist, and so $G$ is non-orderable.
\halmos \end{proof}

Next we'll show:
\begin{Theorem}\label{weeks_circle}
  The fundamental group of the Weeks manifold has no faithful action
  on $S^1$.
\end{Theorem}

\begin{proof}
  Suppose that $G$ acts faithfully on a circle.  Since $H_1(M, \Z/2) =
  0$, the action is orientation preserving.  The Euler class $e$ of
  the action must be nontrivial as the previous theorem shows that $G$
  can't action faithfully on $\R$.  Now $H^2(M) \iso H_1(M) = \Z/5
  \oplus \Z/5$, and so $e$ has order $5$.  By
  Section~\ref{multisections}, this means that $G$ acts faithfully on
  the union of $5$ lines.  Let $N$ be the stabilizer of one of the
  lines. The subgroup $N$ acts faithfully on that line and so is
  left-orderable.  The action on the set of lines is transitive, since
  $5$ is the smallest multiple of $e$ which vanishes; thus $N$ has
  index $5$ in $G$.  To prove the theorem, we'll show that no subgroup
  of $G$ of index $5$ is left-orderable.
  
  One can check that the only subgroups of $G$ of index $5$ are
  normal.  So $N$ is normal and $G/N \iso \Z/5$.  Thus $N$ is the kernel
  of some homomorphism $G \to \Z/5$, and there are $6$ possibilities
  for $N$.  Consider the automorphisms of $G$ given by
  \[
  \phi \maps a \mapsto b, \ b \mapsto a  \mtext{and}
    \psi \maps a \mapsto aB,  \ b \mapsto a.
  \]
  (To check that $\psi$ induces an automorphism, use the relation
  $BaB^2a^2Ba^2B = 1$ mentioned above.)  Together, $\phi$ and $\psi$
  generate a subgroup of the outer automorphism group isomorphic to
  $D_6$ (in fact, this is the whole isometry group of $W$).  
  
  In any event, it is easy to check that under the action of
  $\pair{\phi, \psi}$ there are two orbits of subgroups of index $5$
  in $G$.  Representatives of these orbits are $N_1$, the kernel of
  the homomorphism sending $a \mapsto 0$ and $b \mapsto 1$, and
  $N_2$, the kernel of $a \mapsto 1$ and $b \mapsto -1$.

  First we'll show that $N_1$ is non-orderable.  Let $P$ be the positive
  cone, and assume that $a \in P$.  The proof is completed by
  considering the following cases.  The needed identities in $G$ can be
  checked in number of ways, e.g.~by using automatic groups,
  multiplying matrices, or for the determined, the relations and a
  chalkboard.
  \begin{list}{}{
  \setlength{\labelwidth}{2.6cm}
  \setlength{\leftmargin}{2.9cm}
  \setlength{\labelsep}{0.3cm}
  }

  \item[\textbf{Case $baB \in P$.}\hfil] Then 
    $a (baB)^3 ( a^3 ( baB))^2  a^2 ( baB) (a^3 (baB))^2(baB)^2 = 1 $ in P.
    
  \item[\textbf{Case $bAB \in P$.}\hfil] \ \ 
    \begin{list}{}{
        \setlength{\labelwidth}{3cm}
        \setlength{\leftmargin}{1.2cm}
        \setlength{\labelsep}{0.3cm}
        \setlength{\itemindent}{0cm}
        }
   \item[\textbf{Subcase $abaB \in P$.}\hfil] Then
    \begin{equation*}
       a^2 (abaB)  a^2  (bAB)^2  a  (abaB)  a  ( a  (abaB))^2 
      (a ^2 ( abaB) )^2   a  (abaB) = 1 \  \mbox{in $P$.}
    \end{equation*}  
  
  \item[\textbf{Subcase $bABA \in P$.}]  Set $W = ((bAB)(bABA)(bAB))^2$.  Then 
    \[
    a^2 W (bABA)^2a(bAB)(bABA)(bAB)^2(bABA)W(bAB) = 1\ \mbox{in $P$.}
    \]
    \end{list}
 
\end{list}
This completes the proof that $N_1$ is non-orderable.

Now consider $N_2$.  We can assume that $ab \in P$.  The cases are:
  \begin{list}{}{
  \setlength{\labelwidth}{2.6cm}
  \setlength{\leftmargin}{2.9cm}
  \setlength{\labelsep}{0.3cm}
  }
  \item[\textbf{Case $ba \in P$.}\hfil] Then 
  \begin{list}{}{
  \setlength{\labelwidth}{3cm}
  \setlength{\leftmargin}{1.2cm}
  \setlength{\labelsep}{0.3cm}
  \setlength{\itemindent}{0cm}
  }
   \item[\textbf{Subcase $a^2 b A \in P$.}] Then
   \[
   (ab)  (a^2bA)^2  (ab)^2   (a^2bA)^2  ((ab)^2  (ba)^2)^2  = 1 \  \mbox{in $P$.}
   \]
  
  \item[\textbf{Subcase $aBA^2 \in P$.}]  Then
    \[
(aBA^2)^2(ba)^2(ab)^2(ba)^2(ab)(aBA^2)(ba)^2 
((ab)^2(ba)^2(ab)^2(ba))^2 = 1\ \mbox{in $P$.}
    \]
    \end{list}
    
  \item[ \textbf{Case $AB \in P$.}\hfil] \ 
  \begin{list}{}{
  \setlength{\labelwidth}{3cm}
  \setlength{\leftmargin}{1.2cm}
  \setlength{\labelsep}{0.3cm}
  \setlength{\itemindent}{0cm}
  }
   \item[\textbf{Subcase $a^2bA \in P$.}] Then
     \[
       ((a^2bA)(ab)^2(a^2bA)^2(ab)^2)^2(a^2bA)^2(AB)(ab)  (a^2bA)^2(ab)^2(a^2bA)^2(AB)^2 = 1 \  \mbox{in $P$.}
     \]
  
  \item[\textbf{Subcase $aBA^2 \in P$.}]   Then $(ab)(AB)^3(aBA^2)^3  = 1\ \mbox{in $P$.}$
    \end{list}
\end{list}
This proves that $N_2$ is non-orderable.  This completes the proof that
$G$ does not act faithfully on $S^1$.
\halmos \end{proof}

In fact, one can show more:
\begin{Theorem} 
  Let $G$ be the fundamental group of the Weeks manifold.  Then any
  homomorphism $G \to \homeo(S^1)$ has image which is either trivial or
  $\Z/5\Z$.
\end{Theorem}
\begin{proof}
  Consider a nontrivial action of $G$ on the circle.  Let $K$ be the
  kernel of $G \to \homeo(S^1)$.  We can apply the proofs of
  Theorems~\ref{weeks_nonorder} and \ref{weeks_circle} to $G/K$
  \emph{unless} one of the elements we added to $P$ is in $K$.  In
  other words, we have a contradiction unless one of
  \[
  a, b, aB, bA, B, baB, bAB, abaB, bABA, ba, a^2bA, aBA^2, AB, a^2bA,
  aBA^2
  \]
  is in $K$.  But it is easy to check that the quotient of $G$ by the
  normal closure of any of the above words is $\Z/5\Z$.  So $G/K =
  \Z/5\Z$ and we're done.
\halmos \end{proof}
  
As a corollary of these non-existence results, we get the following:
\begin{Corollary}\label{weeks_nonexistence}
  The Weeks manifold does not admit a taut foliation, a tight essential
  lamination, or a pseudo-Anosov flow.
\end{Corollary}

\begin{proof}
  As the fundamental group of the Weeks manifold can't act on $S^1$,
  the first and last assertions follow immediately from
  Theorem~\ref{universal_circle_action_faithful} and
  Corollary~\ref{cor_anosov_flow}.  For tight essential laminations,
  though, to apply Theorem~\ref{main_thm} we need the additional
  hypothesis that the lamination has solid torus guts.  So it remains
  to show:
  \begin{Claim}
    Let $\Lambda$ be a tight genuine lamination of the Weeks manifold
    $W$.  Then the complementary regions of $\Lambda$ have a
    decomposition where all the gut regions are solid tori.
  \end{Claim}
  
  This will follow from Agol's work on volumes of 3-manifolds with
  tight laminations.  A decomposition of $W - \Lambda$ into
  interstices and guts has \emph{minimal guts} if in every
  complementary region $C$, the interstitial regions are the whole
  characteristic $I$-bundle of $C$.  Such a decomposition is unique up
  to isotopy.  In \cite{AgolVolume}, Agol shows:
  \begin{Theorem}[Agol]\label{volume_bound}
    Let $M$ be a hyperbolic 3-manifold, and $\Lambda$ a tight genuine
    lamination.  Let $G$ be the guts of the minimal decomposition of
    $M - \Lambda$.  Then
    \[
    \vol(M) \ge  -2v_3 \chi(G),
    \]
    where $v_3$ is the volume of the regular ideal tetrahedron in
    $\H^3$.
  \end{Theorem}
  So now let $G_i$ be a gut region of the minimal decomposition of $W -
  \Lambda$.  By Agol's theorem, we must have $\chi(G_i) = 0$ as
  $\vol(W) \approx 0.942 < 2 v_3 \approx 2.02988$.  Thus the boundary
  of $G_i$ consists of tori.  The boundary of $G_i$ is incompressible
  outward into $W$, and so as $W$ is atoroidal, the boundary of $G_i$
  is compressible into $G_i$.  Therefore $G_i$ is a solid torus.  This
  proves the claim, and thus the corollary.
\halmos \end{proof}

\section{Further examples}\label{census_exs}

For the smallest manifolds in the Hodgson-Weeks census \cite{SnapPea}
of closed hyperbolic 3-manifolds we tried to determine which act
faithfully on $\R$.  We looked at the census manifolds of volume $< 3$
which are $\Z/2$-homology spheres (so any action would be orientation
preserving).  There are $128$ such manifolds.  Using the algorithm of
Section~\ref{algorithmic_issues}, we showed that at least $44$ of them
can't act faithfully on $\R$.  Conversely, we  showed
that at least $3$ of them have such actions.  We also found that at
least $60$ have essential laminations.  See the table below for
details.  We would have liked to give more examples where the fundamental
group can't act faithfully on $S^1$, but the only manifold where we
were able to succeed at this was the Weeks manifold.

The algorithm was implemented starting from group presentations
generated by \SnapPea\ and using the automatic groups program KBMAG
\cite{KBMAG} to solve the word problem.

To show manifolds contain essential laminations, we used two techniques.
The first is that for many simple knots, such as two-bridge knots or most
knots under 11 crossings, every nontrivial Dehn surgery contains an
essential lamination (see \cite[pg.~8]{Gabai1997prob} and the references
therein).  In particular, this is true for the knots that appear
in the table below.    In some cases, such as for two-bridge knots, one can
take the lamination to be a taut foliation \cite{Delman95}.

The second technique also uses Dehn filling, but is more complicated
to explain.  Let $M$ be a 3-manifold with torus boundary whose
interior is hyperbolic.  Gabai and Mosher \cite{MosherLams} proved
that $M$ contains a pair of laminations $\Lambda^\pm$ coming from a
pseudo-Anosov flow associated to a finite depth foliation.  When $M$
fibers over $S^1$, these are just the suspensions of the stable and
unstable laminations of the surface homeomorphism.  There is a
\emph{degeneracy slope} in $\bdry M$ associated to $\Lambda^\pm$ with
the following property: for any slope $\alpha$ with $\Delta(\delta,
\alpha) > 1$, the laminations $\Lambda^\pm$ remain essential in the
Dehn filling $M(\alpha)$ (for more see \cite[Thm.~B]{MosherLams} or
the summary of this theory in \cite{Brittenham98}).

In general, it is quite difficult to determine the degeneracy slope.
Here, we used the following trick.  Suppose $M$ has Dehn fillings
$M(\alpha_1),\dots,M(\alpha_n)$ which have finite fundamental group.
The $M(\alpha_i)$ don't contain any essential laminations.  Now
consider some other filling $M(\beta)$.  If $M(\beta)$ also has no
essential lamination, we must have $\Delta(\alpha_i, \delta) \leq 1$
and $\Delta(\beta, \delta) \leq 1$ where $\delta$ is the degeneracy
slope.  So if there does not exist such a $\delta$, we can conclude
that $M(\beta)$ has an essential lamination.  What's more, that
essential lamination in $M(\beta)$ is very full, because the
laminations $\Lambda^\pm$ are themselves very full. 

Below is the table summarizing our findings.  In the Ord column, N
means that the fundamental group is non-orderable, O means that it is
orderable, and blank means unknown.  The Lam column lists an L if the
manifold is known to contain an essential lamination.  If the manifold
is known to contain an essential lamination, the final column gives
the reason.  In that column, $K(p/q)$ means the complement of the
$p/q$ two-bridge knot, and numbers of the form $8_{20}$ refer to the
complements of the corresponding knots in the standard table
\cite{Rolfsen76}.  For the trick with the degeneracy slope, we give the
particular expression as a Dehn filling that was used.

Finally, the reason that the indicated manifolds have orderable
fundamental group is that they have taut foliations and are integral
homology spheres.  In such cases, the action on the universal circle
lifts to a faithful action on $\R$.

{\footnotesize
\begin{longtable}{llcccl}

Name & Volume & Hom & Ord & Lam & Reason for knowing laminar\\
\hline
\endhead
\input small_vol_table
\end{longtable}
}

\newcommand{\etalchar}[1]{$^{#1}$}

\end{document}

%% file: small_vol_table.tex
$m003(-3,1)$ & 0.9427073628 & $\Z/5 + \Z/5$ & N &  & \\
$m003(-2,3)$ & 0.9813688289 & $\Z/5$ &  & L &   Is $m004(5, 1)$ and $m004$ is $K(2/5)$.\\
$m003(-4,3)$ & 1.2637092387 & $\Z/5 + \Z/5$ & N & L &  Degeneracy test as $m003(-4, 3)$.\\
$m004(1,2)$ & 1.3985088842 & $0$ & O  & L &   Is $m004(1, 2)$ and $m004$ is $K(2/5)$.\\
$m003(-4,1)$ & 1.4236119003 & $\Z/35$ & N &  & \\
\hline
$m004(3,2)$ & 1.4406990067 & $\Z/3$ &  & L &   Is $m004(3, 2)$ and $m004$ is $K(2/5)$.\\
$m004(7,1)$ & 1.4637766449 & $\Z/7$ &  & L &   Is $m004(7, 1)$ and $m004$ is $K(2/5)$.\\
$m004(5,2)$ & 1.5294773294 & $\Z/5$ &  & L &   Is $m004(5, 2)$ and $m004$ is $K(2/5)$.\\
$m003(-5,3)$ & 1.5435689115 & $\Z/35$ & N & L &  Degeneracy test as $m003(-5, 3)$.\\
$m007(1,2)$ & 1.5435689115 & $\Z/21$ & N & L &  Degeneracy test as $m011(3, 2)$.\\
\hline
$m007(4,1)$ & 1.5831666606 & $\Z/21$ & N &  & \\
$m007(3,2)$ & 1.5831666606 & $\Z/3 + \Z/9$ & N &  & \\
$m006(-3,2)$ & 1.6496097158 & $\Z/15$ & N & L &  Degeneracy test as $m006(-3, 2)$.\\
$m015(5,1)$ & 1.7571260292 & $\Z/7$ &  & L &   Is $m015(5, 1)$ and $m015$ is $K(-2/7)$.\\
$m007(-3,2)$ & 1.8243443222 & $\Z/3 + \Z/3$ & N & L &  Degeneracy test as $m007(-3, 2)$.\\
\hline
$m016(-3,2)$ & 1.8854147256 & $\Z/39$ & N & L &  Degeneracy test as $m016(-3, 2)$.\\
$m017(-3,2)$ & 1.8854147256 & $\Z/7 + \Z/7$ & N & L &  Degeneracy test as $m017(-3, 2)$.\\
$m006(3,2)$ & 1.8859142560 & $\Z/45$ & N & L &  Degeneracy test as $m006(3, 2)$.\\
$m011(2,3)$ & 1.9122102501 & $0$ &  & L &   Is $m222(-2, 1)$ and $m222$ is $8_{20}$.\\
$m006(4,1)$ & 1.9222971095 & $\Z/35$ & N &  & \\
\hline
$m006(-2,3)$ & 1.9537083154 & $\Z/35$ & N & L &  Degeneracy test as $m006(-2, 3)$.\\
$m006(2,3)$ & 1.9627376578 & $\Z/55$ & N & L &  Degeneracy test as $m006(2, 3)$.\\
$m017(-1,3)$ & 1.9627376578 & $\Z/7 + \Z/7$ & N &  & \\
$m023(-4,1)$ & 2.0143365838 & $\Z/3 + \Z/3$ &  &  & \\
$m007(5,2)$ & 2.0259452819 & $\Z/33$ & N &  & \\
\hline
$m006(-5,2)$ & 2.0288530915 & $\Z/5$ &  & L &   Is $m015(1, 2)$ and $m015$ is $K(-2/7)$.\\
$m036(-3,2)$ & 2.0298832128 & $\Z/3 + \Z/15$ &  &  & \\
$m007(-6,1)$ & 2.0555467489 & $\Z/3 + \Z/3$ &  &  & \\
$m007(-5,2)$ & 2.0656708385 & $\Z/3$ &  & L &   Is $m015(-1, 2)$ and $m015$ is $K(-2/7)$.\\
$m015(-5,1)$ & 2.1030952907 & $\Z/3$ &  & L &   Is $m015(-5, 1)$ and $m015$ is $K(-2/7)$.\\
\hline
$m016(3,2)$ & 2.1145676931 & $\Z/33$ & N & L &  Degeneracy test as $m016(3, 2)$.\\
$m015(3,2)$ & 2.1145676931 & $\Z/7$ &  & L &   Is $m015(3, 2)$ and $m015$ is $K(-2/7)$.\\
$m011(4,1)$ & 2.1243017573 & $\Z/43$ &  &  & \\
$m017(1,3)$ & 2.1557385676 & $\Z/35$ & N &  & \\
$m011(-2,3)$ & 2.1557385676 & $\Z/53$ & N & L &  Degeneracy test as $m011(-2, 3)$.\\
\hline
$m034(4,1)$ & 2.1847555751 & $\Z/7$ &  & L &   Is $s385(-2, 1)$ and $s385$ is $10_{125}$.\\
$m034(-4,1)$ & 2.1959641187 & $\Z/25$ & N &  & \\
$m011(-3,2)$ & 2.2082823597 & $\Z/57$ & N & L &  Degeneracy test as $m011(-3, 2)$.\\
$m011(4,3)$ & 2.2102443409 & $\Z/25$ &  & L &  Degeneracy test as $m011(4, 3)$.\\
$m011(1,4)$ & 2.2109517391 & $\Z/23$ &  & L &  Degeneracy test as $m011(1, 4)$.\\
\hline
$m015(-3,2)$ & 2.2267179039 & $0$ & O & L &   Is $m015(-3, 2)$ and $m015$ is $K(-2/7)$.\\
$m015(7,1)$ & 2.2267179039 & $\Z/9$ &  & L &   Is $m015(7, 1)$ and $m015$ is $K(-2/7)$.\\
$m038(1,2)$ & 2.2597671326 & $0$ &  & L &   Is $m372(-2, 1)$ and $m372$ is $9_{46}$.\\
$m015(5,2)$ & 2.2662435733 & $\Z/9$ &  & L &   Is $m015(5, 2)$ and $m015$ is $K(-2/7)$.\\
$m026(-4,1)$ & 2.2726318636 & $\Z/13$ &  &  & \\
\hline
$m011(-1,4)$ & 2.2757758101 & $\Z/49$ & N & L &  Degeneracy test as $m011(-1, 4)$.\\
$m023(-3,2)$ & 2.2944383001 & $\Z/3$ &  & L &   Is $m032(5, 1)$ and $m032$ is $K(-2/9)$.\\
$m038(-5,1)$ & 2.3126354033 & $\Z/17$ &  &  & \\
$m017(-5,1)$ & 2.3188118677 & $\Z/7 + \Z/7$ & N & L &  Degeneracy test as $m022(-3, 2)$.\\
$m016(-5,1)$ & 2.3188118677 & $\Z/23$ &  &  & \\
\hline
$m019(4,1)$ & 2.3207602675 & $\Z/7$ &  & L &   Is $m199(3, 1)$ and $m199$ is $9_{42}$.\\
$m022(1,3)$ & 2.3380401178 & $\Z/35$ & N &  & \\
$m016(-1,4)$ & 2.3522069054 & $\Z/73$ & N & L &  Degeneracy test as $m026(2, 3)$.\\
$m017(-1,4)$ & 2.3522069054 & $\Z/63$ & N &  & \\
$m019(-2,3)$ & 2.3641969332 & $\Z/63$ & N & L &  Degeneracy test as $m019(-2, 3)$.\\
\hline
$m022(5,1)$ & 2.3705924006 & $\Z/3 + \Z/7$ &  &  & \\
$m019(-4,1)$ & 2.3803358221 & $\Z/41$ & N & L &  Degeneracy test as $m026(-2, 3)$.\\
$m022(5,2)$ & 2.4224625169 & $\Z/7$ &  & L &   Is $m032(-5, 1)$ and $m032$ is $K(-2/9)$.\\
$m019(4,3)$ & 2.4444077795 & $\Z/27$ &  & L &  Degeneracy test as $m019(4, 3)$.\\
$m022(-1,3)$ & 2.4540294422 & $\Z/7 + \Z/7$ &  &  & \\
\hline
$m026(4,1)$ & 2.4631393944 & $\Z/51$ &  &  & \\
$m029(-3,2)$ & 2.4682321967 & $\Z/5 + \Z/9$ & N &  & \\
$m036(3,2)$ & 2.4682321967 & $\Z/3 + \Z/9$ & N & L &  Degeneracy test as $m036(3, 2)$.\\
$m022(-5,1)$ & 2.4878225918 & $\Z/7 + \Z/7$ & N &  & \\
$m023(-6,1)$ & 2.4903791858 & $\Z/15$ &  &  & \\
\hline
$m038(3,2)$ & 2.5026593054 & $\Z/5$ &  & L &   Is $m289(2, 1)$ and $m289$ is $K(-3/11)$.\\
$m034(-5,1)$ & 2.5065758445 & $\Z/29$ & N &  & \\
$m034(5,1)$ & 2.5144043349 & $\Z/11$ &  &  & \\
$m070(-3,1)$ & 2.5274184773 & $\Z/11$ &  & L &  Degeneracy test as $m117(-3, 2)$.\\
$m038(-5,2)$ & 2.5274184773 & $\Z/19$ &  &  & \\
\hline
$m036(-5,1)$ & 2.5274184773 & $\Z/33$ &  &  & \\
$m030(5,2)$ & 2.5303032876 & $\Z/63$ & N &  & \\
$m023(-5,2)$ & 2.5415850101 & $\Z/3 + \Z/3$ &  &  & \\
$m038(5,1)$ & 2.5495466001 & $\Z/13$ &  &  & \\
$m026(-5,1)$ & 2.5667347900 & $\Z/21$ &  &  & \\
\hline
$m160(1,2)$ & 2.5689706009 & $\Z/3 + \Z/5$ &  &  & \\
$m036(-1,3)$ & 2.5751620736 & $\Z/57$ &  &  & \\
$m030(1,3)$ & 2.5854830480 & $\Z/7 + \Z/7$ &  & L &  Degeneracy test as $m030(1, 3)$.\\
$m160(-3,2)$ & 2.5953875937 & $\Z/3 + \Z/9$ & N &  & \\
$m036(-5,2)$ & 2.6095439552 & $\Z/51$ &  &  & \\
\hline
$m027(-4,1)$ & 2.6122234482 & $\Z/77$ &  &  & \\
$m027(4,3)$ & 2.6172815707 & $\Z/25$ &  & L &  Degeneracy test as $m027(4, 3)$.\\
$m081(1,3)$ & 2.6244624283 & $\Z/37$ & N &  & \\
$m036(5,1)$ & 2.6285738915 & $\Z/3$ &  & L &   Is $s580(-2, 1)$ and $s580$ is $10_{145}$.\\
$m032(5,2)$ & 2.6294053953 & $0$ & O  & L &   Is $m032(5, 2)$ and $m032$ is $K(-2/9)$.\\
\hline
$m034(-1,3)$ & 2.6414714456 & $\Z/31$ &  & L &  Degeneracy test as $m034(-1, 3)$.\\
$m036(1,3)$ & 2.6536080625 & $\Z/51$ &  & L &  Degeneracy test as $m082(-3, 2)$.\\
$m034(-2,3)$ & 2.6555425236 & $\Z/35$ &  &  & \\
$m034(1,3)$ & 2.6646126469 & $\Z/23$ &  & L &  Degeneracy test as $m034(1, 3)$.\\
$m160(2,1)$ & 2.6735274161 & $\Z/3$ &  & L &   Is $m372(2, 1)$ and $m372$ is $9_{46}$.\\
\hline
$m032(7,1)$ & 2.6822267321 & $\Z/5$ &  & L &   Is $m032(7, 1)$ and $m032$ is $K(-2/9)$.\\
$m069(4,1)$ & 2.6954841673 & $\Z/65$ & N & L &  Degeneracy test as $m081(-3, 2)$.\\
$m069(-1,3)$ & 2.6954841673 & $\Z/39$ &  &  & \\
$m030(5,3)$ & 2.7067833105 & $\Z/77$ & N & L & Is Haken.  See \cite{Dunfield:haken}.\\
$m120(-3,2)$ & 2.7124588084 & $0$ &  & L &   Is $m199(-3, 1)$ and $m199$ is $9_{42}$.\\
\hline
$m116(-1,3)$ & 2.7589634387 & $\Z/7$ &  & L &   Is $s580(2, 1)$ and $s580$ is $10_{145}$.\\
$m081(-1,3)$ & 2.7725163132 & $\Z/59$ & N &  & \\
$m160(-2,3)$ & 2.8022537823 & $\Z/3 + \Z/11$ &  &  & \\
$m221(3,1)$ & 2.8281220883 & $\Z/21$ &  &  & \\
$m142(3,2)$ & 2.8281220883 & $\Z/19$ &  &  & \\
\hline
$m206(1,2)$ & 2.8281220883 & $\Z/5$ &  &  & \\
$m082(2,3)$ & 2.8458961160 & $\Z/83$ & N &  & \\
$m070(4,3)$ & 2.8472238006 & $\Z/85$ & N &  & \\
$m069(4,3)$ & 2.8472238006 & $\Z/99$ &  &  & \\
$m137(-5,1)$ & 2.8656302333 & $0$ &  &  & \\
\hline
$m070(-2,3)$ & 2.8669017766 & $\Z/61$ & N & L &  Degeneracy test as $m070(-2, 3)$.\\
$m069(-2,3)$ & 2.8669017766 & $\Z/27$ &  & L &  Degeneracy test as $m069(-2, 3)$.\\
$m069(-4,1)$ & 2.8733431176 & $\Z/31$ &  &  & \\
$m070(-4,1)$ & 2.8733431176 & $\Z/7$ &  &  & \\
$m100(2,3)$ & 2.8824943873 & $\Z/85$ &  &  & Is Haken. See \cite{Dunfield:haken}. \\
\hline
$m082(-2,3)$ & 2.9027039980 & $\Z/79$ &  & L &  Degeneracy test as $m082(-2, 3)$.\\
$m221(-1,2)$ & 2.9133321143 & $\Z/7$ &  &  & \\
$m116(1,3)$ & 2.9169341134 & $\Z/41$ &  &  & \\
$m120(-5,1)$ & 2.9356518985 & $\Z/17$ &  &  & \\
$m078(2,3)$ & 2.9398104423 & $\Z/37$ &  &  & \\
\hline
$m145(2,3)$ & 2.9400386172 & $\Z/47$ & N &  & \\
$m078(5,2)$ & 2.9438596478 & $\Z/43$ &  &  & \\
$m249(3,1)$ & 2.9545326040 & $\Z/3 + \Z/5$ &  &  & \\
$m145(3,2)$ & 2.9582502906 & $\Z/13$ &  &  & \\
$m117(3,2)$ & 2.9605565159 & $\Z/53$ & N &  & \\
\hline
$m117(-5,1)$ & 2.9607151670 & $\Z/19$ &  &  & \\
$m154(2,3)$ & 2.9670703390 & $\Z/77$ &  &  & \\
$m078(-2,3)$ & 2.9696321386 & $\Z/17$ &  &  & \\
$m100(-2,3)$ & 2.9709840073 & $\Z/77$ &  & L &  Degeneracy test as $m100(-2, 3)$.\\
$m117(1,3)$ & 2.9760925194 & $\Z/55$ &  &  & \\
\hline
$m078(-5,2)$ & 2.9769925267 & $\Z/7$ &  & L &   Is $m199(-1, 2)$ and $m199$ is $9_{42}$.\\
$m159(3,2)$ & 2.9781624873 & $\Z/35$ &  &  & \\
$m137(5,1)$ & 2.9868370451 & $0$ &  &  & \\